\documentclass[11pt,english,leqno]{article}
\usepackage[a4paper]{geometry}
\usepackage{setspace}

\RequirePackage{amsthm,amsmath,amsfonts,amssymb}
\RequirePackage[authoryear]{natbib}
\RequirePackage{color}
\RequirePackage{mathtools}
\RequirePackage{caption}

\theoremstyle{plain}
\newtheorem{thm}{Theorem}
\newtheorem{lem}{Lemma}

\newtheorem{assumption}{Assumption}
\newtheorem{myprocedure}{Procedure}
\newtheorem*{myassumption1*}{Assumption GAN}
\newtheorem*{myassumption2*}{Assumption Smooth GAN}
\newtheorem*{myassumption3*}{Assumption Minimax}
\newtheorem*{myassumption4*}{Assumption Smooth Minimax}

\newenvironment{lyxlist}[1]
	{\begin{list}{}
		{\settowidth{\labelwidth}{#1}
		 \setlength{\leftmargin}{\labelwidth}
		 \addtolength{\leftmargin}{\labelsep}
		 }}
	{\end{list}}

\DeclareMathOperator{\CS}{CS}
\DeclareMathOperator{\Var}{Var}

\newcommand{\smalluparrow}{%
  \mathrel{\nonscript\mkern-1.2mu\mkern1.2mu{\scriptscriptstyle\uparrow}}%
}

\begin{document}
\title{\vspace{20pt}
Statistical inference for generative adversarial networks\\
and other minimax problems\,%
\thanks{%
First version available as arXiv:2104.10601 in April 2021.
The author thanks
the Research Council of Finland,
Foundation for the Advancement of Finnish Securities Markets,
and OP Group Research Foundation
for financial support,
CSC -- IT Center for Science, Finland, for computational resources,
and the anonymous referees for useful comments and suggestions.
Contact address: Mika Meitz, Department of Economics,
University of Helsinki, P. O. Box 17, FI\textendash 00014 University
of Helsinki, Finland; e-mail: mika.meitz@helsinki.fi.%
}\vspace{20pt}
}
\author{Mika Meitz\\\small{University of Helsinki} \vspace{20pt}
}
\date{December 2023}
\maketitle
\begin{abstract}
\noindent%
This paper studies generative adversarial networks (GANs)
from the perspective of statistical inference.
A GAN is a popular machine learning method in which the parameters of two neural networks, 
a generator and a discriminator, are estimated to solve a particular minimax problem.
This minimax problem typically has a multitude of solutions 
and the focus of this paper are the statistical properties of these solutions.
We address two key statistical issues for the generator and discriminator network parameters, 
consistent estimation and confidence sets.
We first show that the set of solutions to the sample GAN problem is a 
(Hausdorff) consistent estimator of the set of solutions to the corresponding population GAN problem.
We then devise a computationally intensive procedure to form confidence sets 
and show that these sets contain the population GAN solutions 
with the desired coverage probability. 
Small numerical experiments and a Monte Carlo study
illustrate our results and verify our theoretical findings.
We also show that our results apply in general minimax problems
that may be non-convex, non-concave, and have multiple solutions.

\bigskip{}
\bigskip{}
\bigskip{}

\noindent\textbf{MSC2020 subject classifications:} Primary 62F12; secondary 68T07.

\bigskip{}


\end{abstract}
\vfill{}

\pagebreak{}

\section{Introduction}

A generative adversarial network (GAN) is a machine learning method
introduced in \citet{goodfellow2014generative}. The basic purpose
of a GAN is to learn how to generate synthetic data based on a training
set of real-world examples. While a traditional use of GANs has been
to generate authentic-looking photographs based on real example images
(for a recent example, see \citealp{karras2020style}), they are now
in use across various scientific disciplines. To give brief examples,
variants of GANs have been used 
in biology to generate synthetic protein-coding DNA sequences (\citealp{gupta2019feedback}), 
in physics to simulate subatomic particle collisions at the Large Hadron Collider (\citealp{paganini2018accelerating}),
in astronomy to de-noise images of galaxies (\citealp{schawinski2017generative}),
and in medicine to improve near-infrared fluorescence imaging (\citealp{mae2021deep}). 
GANs have also been used in popular culture and arts 
to re-create video games (\citealp{kim2020learning}),
to make computer-generated art (\citealp{miller2019artist}), 
and to compose music (\citealp{briot2020deep}). 
To put the recent popularity of GANs into perspective,
the original \citet{goodfellow2014generative} article
received over 10.000 Google Scholar citations in the year 2020 alone.
For recent surveys of GANs and for further references, see
\citet{creswell2018generative}, \citet{pan2019recent}, or \citet{goodfellow2020generative}.

A GAN typically comprises two neural networks called a generator
and a discriminator. (For background material on neural networks,
see the books of \citet{bishop2006pattern} or \citet{goodfellow2016deep}.)
The generator network produces synthetic data whose distribution aims
to mimic that of real data, while the discriminator network evaluates
whether the data produced by the generator is real or fake. Suppose
the observed vector $x$ mathematically represents an observed image,
DNA sequence, or some other object of interest. This $x$ is viewed
as a realization from an underlying random vector $X$ whose distribution
remains unknown to us. Obtaining a large number of realizations from
$X$ may be difficult or costly, and the researcher desires to produce
synthetic replicas in an easy manner. These replicas are the output
of a generating mechanism $G(z,\gamma)$ that in practice is a complicated
neural network. This generator takes as an input noise variables $z$
that are drawn from some underlying random vector $Z$ and depends
on parameters $\gamma$ that are tuned so that the output $G(z,\gamma)$
would be a close replica of $x$. To assess the quality of the synthetic
data produced, a discriminating mechanism $D(\cdot,\delta)$ indicates
how likely it is an input, whether original data $x$ or a replica
$G(z,\gamma)$, is real data from the underlying distribution $X$.
In practice, the discriminator $D(\cdot,\delta)$ is again a complicated
neural network with parameters $\delta$ to be estimated. The `adversarial
principle' of generative adversarial networks works as follows. Given
an original $x$ and any synthetic $G(z,\gamma)$, the discriminator
aims to give a high rating $D(x,\delta)$ to the real $x$ and a low
rating $D(G(z,\gamma),\delta)$ to the artificial $G(z,\gamma)$ by
choosing $\delta$ to maximize the objective function
\[
D( x,\delta)
\times
(\thinspace1-D( G(z,\gamma),\delta)\thinspace)
\]
for any given fixed $\gamma$ at a time. In contrast, the generator
aims to make $G(z,\gamma)$ as hard to distinguish from $x$ as possible
by choosing the $\gamma$ to minimize the discriminator's maximized
objective function. 

To formalize the above discussion, the GAN problem can be expressed
as the minimax problem 
\begin{equation}
\adjustlimits\inf_{\gamma\in\Gamma}\sup_{\delta\in\Delta}f(\gamma,\delta)
\text{\ \ with\ \ }
f(\gamma,\delta)=\mathbb{E}[\ln(D(X,\delta))]+\mathbb{E}[\ln(1-D(G(Z,\gamma),\delta))],
\label{eq:IntroGANpopula}
\end{equation}
where the infimum and supremum are taken over sets $\Gamma$ and $\Delta$
denoting the ranges of permissible values for $\gamma$ and $\delta$,
respectively, and $\mathbb{E}$ denotes expectation with respect to
the joint distribution of $X$ and $Z$. The formalities will be discussed
in detail in Section 2. For now it suffices to remark that $X$ and
$Z$ are typically of rather large dimension, that the neural networks
$G(\cdot,\cdot)$ and $D(\cdot,\cdot)$ are in practice quite complicated
and parameterized via very high-dimensional $\gamma$ and $\delta$,
and that in GAN applications $f$ is as a rule non-convex and non-concave
(in contrast to the traditional convex-concave minimax setting in
which $f(\cdot,\delta)$ is convex for all fixed $\delta\in\Delta$
and $f(\gamma,\cdot)$ is concave for all fixed $\gamma\in\Gamma$).
Given solutions, say $\gamma_{0}$ and $\delta_{0}$, to problem (\ref{eq:IntroGANpopula}),
especially $\gamma_{0}$ is of interest as $G(Z,\gamma_{0})$ gives
the researcher a mechanism to produce the desired synthetic data.
This description of GANs corresponds to the original formulation in
\citet{goodfellow2014generative}. The original GAN and its numerous
extensions and variants have in recent years attracted remarkable
interest in applications in which having access to large quantities of
synthetic data is beneficial. Examples of such applications were mentioned
above, and the surveys listed above contain further references as
well as details of the many extensions of GANs.

The main object of interest in this paper is the GAN minimax problem (\ref{eq:IntroGANpopula}) and its solutions. 
On a more general level, recent machine learning literature has experienced 
a notable surge of interest in general minimax problems. As GANs have arguably been one of the major reasons for this, we focus on the GAN case as a prominent representative example of minimax problems; 
other minimax problems will be discussed in Section 6.

In practical GAN applications, one would use a large training sample
of observations to estimate the GAN. Let the $n$ observations correspond
to independent and identically distributed (IID) random vectors $X_{1},\ldots,X_{n}$,
all distributed as the variable $X$, and suppose the $n$\, IID random
vectors $Z_{1},\ldots,Z_{n}$ have the same distribution as the noise
variable $Z$. Then the objective is to solve the sample minimax problem
(cf.~\citealp{goodfellow2014generative}; \citealp{biau2020some})
\begin{equation}
\adjustlimits\inf_{\gamma\in\Gamma}\sup_{\delta\in\Delta}\hat{f}_{n}(\gamma,\delta)
\text{\ \ with\ \ }
\hat{f}_{n}(\gamma,\delta)=\frac{1}{n}\sum_{i=1}^{n}\ln(D(X_{i},\delta))
+ \frac{1}{n}\sum_{i=1}^{n}\ln(1-D(G(Z_{i},\gamma),\delta)).
\label{eq:IntroGANsample}
\end{equation}
The primary objective in GAN applications in the machine learning
literature is to find solutions to this sample problem, or to `train'
the GAN. In statistical terminology, this corresponds to estimating
the GAN parameters $\gamma$ and $\delta$. (Of course, choosing the
parametric forms of the neural networks $G(\cdot,\cdot)$ and $D(\cdot,\cdot)$,
or specifying the network architecture, is done before this.) 
Devising algorithms that solve the sample
minimax problem (\ref{eq:IntroGANsample}) is very challenging and
a large body of machine learning literature on GANs focuses on this.
Discussions of convergence and stability of various GAN training algorithms
can be found for instance in \citet{arjovsky2017towards}, \citet{nagarajan2017gradient},
and \citet{mescheder2018training}; more recent contributions include
\citet{diakonikolas2021efficient}, \citet{fiez2021local}, and \citet{mangoubi2021greedy}.
Another important theoretical aspect of GANs studied in the machine learning
literature is how well the distribution of $G(Z,\gamma)$ can approximate the target distribution of $X$;
contributions to this line of research include 
\citet{arora2017generalization}, \citet{liu2017approximation}, \citet{singh2018nonparametric},
\citet{lu2020universal}, \citet{biau2021some}, \citet{liang2021well}, and \citet{schreuder2021statistical}.

This paper studies 
the GAN minimax problem
from 
the perspective of statistical inference.
We do not address algorithmic issues and simply
assume a method is available for solving the sample GAN problem (\ref{eq:IntroGANsample}).
The focus of this paper are the statistical properties of these solutions,
say $\hat{\gamma}_{n}$ and $\hat{\delta}_{n}$, as estimators of
the solutions $\gamma_{0}$ and $\delta_{0}$ to the population problem (\ref{eq:IntroGANpopula}). 
Such questions of statistical inference are quite
orthogonal to much of the machine learning literature on GANs: 
In most GAN applications the ability to produce synthetic data is the ultimate goal 
and, as discussed above, theoretical works on GANs often focus on convergence of algorithms or the 
ability of GANs to mimic the target distribution.
Our study can be seen as complementary to these existing works.
Although the parameters $\gamma$ and $\delta$ have no particular interpretation
in GANs, studying properties of their estimators nevertheless contributes 
to a more complete understanding of the statistical properties of GANs.
In the broader context of general neural network models,
it is also typically the case that
the network parameters (weights) are not of major interest.
Nevertheless, statistical questions regarding these parameters
have been explored in the literature. 
To give a few examples, research in this direction includes
early frequentist works on neural network parameter estimation
(\citealp{white1989learning}; \citealp{white1989some}), as well as
more recent works on Bayesian posterior distribution of the neural network parameters 
(\citealp{blundell2015weight}; \citealp{izmailov2021what})
and tests to assess the statistical significance of neural network variables
(\citealp{horel2020significance}; \citealp{fallahgoul2024asset}).
The present paper relates to this earlier literature and 
views the GAN minimax solutions as objects of statistical interest.

The study of statistical inference 
for the GAN solutions 
was initiated in the
recent important paper of \citet{biau2020some}.
These authors were the first to consider the consistency and asymptotic
normality of the sample GAN solutions $\hat{\gamma}_{n}$ and $\hat{\delta}_{n}$
as estimators of the population GAN solutions $\gamma_{0}$ and $\delta_{0}$.
A key assumption these authors make is that the population GAN problem
(\ref{eq:IntroGANpopula}) has a single, unique solution $(\gamma_{0},\delta_{0})$.
However, as \citet[pp.~1560--1561]{biau2020some} acknowledge,
 this assumption is unrealistic in practical GAN applications.
On the contrary,
as the objective function $f(\gamma,\delta)$ in (\ref{eq:IntroGANpopula})
is parameterized using complicated neural networks $G$ and $D$,
it is essentially guaranteed to have an extremely large number of
solutions (for related discussion, see \citealp[Sec 8.2]{goodfellow2016deep}).
One reason for this prevalence of multiple solutions is the inherent
non-identifiability of heavily parameterized neural networks. 

In this paper, we consider statistical inference for GANs in the empirically
relevant case of multiple solutions. 
The theoretical framework required for this is rather different from the one employed in \citet{biau2020some}.
We focus on two key issues, consistent
estimation and confidence sets, and the assumptions we make are weak
and hold in many real GAN applications. To briefly describe our results,
let $\Theta_{0}$ and $\hat{\Theta}_{n}$ denote the sets of solutions
to the population minimax problem (\ref{eq:IntroGANpopula}) and the
sample GAN problem (\ref{eq:IntroGANsample}), respectively. 
We first consider an appropriately defined notion of consistency 
for the set-valued estimator $\hat{\Theta}_{n}$, namely Hausdorff consistency
(precise definition will be given in Section 3). Without any restrictions
on the (potentially infinite) number of solutions, we show that $\hat{\Theta}_{n}$
is a Hausdorff consistent estimator of $\Theta_{0}$. We then consider
confidence sets, random sets that contain $\Theta_{0}$
with a prespecified coverage probability. In the traditional point-identified
setting in which $\Theta_{0}$ is a singleton, confidence sets would
often be formed based on the asymptotic distribution of an estimator
of the parameters of interest. In the present set-identified case,
we follow an alternative approach. 
We devise a computationally intensive resampling procedure based on appropriate lower contour sets of a particular criterion function to form the confidence sets,
and show that these sets contain the population GAN solutions 
with (at least) the desired coverage probability. 
This turns out to be technically challenging, and the details will be given in Section 4 and Appendix B.

The theoretical developments in this paper build on two strands of literature. 
On the one hand, our results for consistent estimation and for
constructing confidence sets are based on recent developments in estimation
and inference for general set-identified parameters in the partial
identification literature in econometrics; our approach relies in
particular on the pioneering work of \citet{chernozhukov2007estimation}
and \citet{romano2010inference} (for further related references, see the
recent survey of this literature in \citealp{molinari2020microeconometrics}).
Note that conventional theory of point-identified extremum estimators 
(see, e.g., \citealp[Ch 5]{vanderVaart1998asymptotic})
is not applicable in our setting.
The research in the partial identification literature does not, however,
consider minimax problems such as (\ref{eq:IntroGANpopula}) or GANs, 
and adjusting to the minimax setting requires some work.
On the other hand, 
general minimax problems (but no GANs)
have previously been considered in the stochastic programming literature in operations
research. Our results are particularly closely related to those
in \citet{shapiro2008asymptotics} and \citet[Sec 5.1.4]{shapiro2009lectures}
who consider consistent estimation (but no confidence sets) in minimax
problems, as well as to the foundational results for minimization problems
in \citet{shapiro1991asymptotic}. Further comparison to all these
previous works will be given later in the paper. 

The present paper contributes to the statistical literature in several ways.
Although other theoretical aspects of GANs have been considered previously, 
the present paper provides the first results for statistical inference in GAN minimax problems
in the practically relevant case of multiple solutions. 
Our results also apply to other general minimax problems, not just to GANs.
The confidence sets we provide for the solutions of minimax problems are novel
and, in contrast to previous literature,  
the minimax problems may be non-convex, non-concave, and have multiple solutions
(see Section 6 for details).
Furthermore, as a technical device in our proofs, we establish new Hadamard directional differentiability
results of certain mappings related to GANs and general minimax problems, and these results
may be of independent interest (see Lemma 3 in Appendix B).

Finally, it should be noted that the results in this paper are asymptotic, 
taking the sample size $n$ to infinity while keeping the dimensions of the parameters $\gamma$ and $\delta$ fixed.
This traditional asymptotic framework is appropriate in minimax problems with a moderate number of parameters but is somewhat idealized in GAN applications where the number of parameters may be extremely large or even substantially exceeding $n$.
Our asymptotic results aim to contribute towards a better understanding of 
statistical inference for set-identified parameters in GANs and other minimax problems;
exploring non-asymptotic results would also be very interesting but would require a different mathematical framework and is beyond the scope of this paper
(we note that, to the best of our knowledge, non-asymptotic theory of estimation for set-identified parameters has not yet been developed in the literature).

The rest of the paper is organized as follows. The next section sets
the stage by considering the GAN minimax problem more formally. Consistent
estimation of $\Theta_{0}$ is the topic of Section 3, while Section
4 discusses confidence sets for $\Theta_{0}$. 
These results are illustrated in small numerical experiments in Section 5.
Other general minimax problems are discussed in Section 6, and
Section 7 concludes.
All technical derivations and proofs are relegated to Appendices A--C
(with Appendix B containing the most interesting results).

\section{The GAN minimax problem}

We now consider the GAN problem outlined in the Introduction more
formally. In what follows, all the random quantities are defined on
some appropriate underlying probability space but typically there
is no need to emphasize this. The available data corresponds to a
random vector $X$ taking values in some Euclidean subset $\mathcal{X}\subseteq\mathbb{R}^{d_{X}}$
and whose distribution remains unknown. The noise variables $Z$ with
values in $\mathcal{Z}\subseteq\mathbb{R}^{d_{Z}}$ come from a distribution
chosen by the researcher, often the multivariate uniform or Gaussian
distribution. The generator function $G(\cdot,\cdot)$ is a neural
network that transforms the noise variables $Z$ into synthetic data
$G(Z,\gamma)$. The discriminator function $D(\cdot,\cdot)$ is another 
neural network indicating how likely it is that a real observation
$x$ or a replica $G(z,\gamma)$ comes from the distribution of $X$.
The following assumption summarizes these concepts from a technical
point of view.

\begin{myassumption1*}
\mbox{}
\begin{lyxlist}{(a)}
\item [{(a)}] Suppose $X_{1},\ldots,X_{n}$ and $Z_{1},\ldots,Z_{n}$ are independent and identically distributed random vectors with the same distributions as $X$ and $Z$ and taking values in $\mathcal{X}\subseteq\mathbb{R}^{d_{X}}$ and $\mathcal{Z}\subseteq\mathbb{R}^{d_{Z}}$, respectively.  
\item [{(b)}] The set $\Theta=\Gamma\times\Delta\subseteq\mathbb{R}^{d_{\gamma}+d_{\delta}}$ is compact and non-empty.  
\item [{(c)}] The generator function $G:\mathcal{Z}\times\Gamma\to\mathcal{X}$, the discriminator function $D:\mathcal{X}\times\Delta\to(0,1)$, and the function $F:\mathcal{X}\times\mathcal{Z}\times\Gamma\times\Delta\to\mathbb{R}$ defined by \[ F(x,z,\gamma,\delta)=\ln(D(x,\delta))+\ln(1-D(G(z,\gamma),\delta)) \] are such that $F(X,Z,\gamma,\delta)$ is measurable for all $(\gamma,\delta)\in\Theta$ and continuous on $\Theta$ with probability one.
\end{lyxlist} 
\end{myassumption1*}

This assumption contains some minimal requirements for the original
GAN framework of \citet{goodfellow2014generative}. The IID assumption
in part (a) is standard in the GAN setting, and the same holds for
the compactness assumption for the permissible parameter space in part (b). 
Part (c) is a minimal continuity
assumption. Similar assumptions have been used for instance by \citet{biau2020some}.

Descriptions of neural network architectures that can be used to specify the
generator $G$ and the discriminator $D$ can be found in the book
of \citet{goodfellow2016deep} and in the many references therein. 
Often the generator and
discriminator networks would satisfy additional differentiability and moment condition
assumptions. Although such extra assumptions are not necessarily required for
the results that follow, such a `smooth GAN' setting serves as a convenient
example that will be used to illustrate our results. In the following
assumption and in what follows, we use either notation $\theta$
or notation $(\gamma,\delta)$ for the elements of $\Theta=\Gamma\times\Delta$;
for instance, $F(x,z,\theta)$ and $F(x,z,\gamma,\delta)$ 
are used interchangeably.
Also, $\left|\cdot\right|$ denotes the Euclidean norm.

\begin{myassumption2*}
Assumption GAN holds, 
the function $F(x,z,\theta)$ is continuously differentiable on 
an open convex set $\Theta^{\ast}$ containing
$\Theta$ for all $(x,z)\in\mathcal{X}\times\mathcal{Z}$
with a square integrable derivative
$(\mathbb{E}[
\sup_{\theta\in\Theta^{\ast}}
|\partial F(X,Z,\theta)/\partial\theta|^{2}]<\infty)$,
and\/ $\inf_{\theta\in\Theta}\Var[F(X,Z,\theta)]>0$. 
\end{myassumption2*}

These conditions require the GAN problem to be somewhat more well-behaved. 
Differentiability of the generator and discriminator networks is a common requirement for training methods employed in GAN applications, 
such as variants of the gradient descent-ascent algorithm
(the set $\Theta^{\ast}$ is introduced to ensure derivatives are well-defined
throughout the proofs).
The mild integrability condition facilitates checking later high-level assumptions and holds, for example, 
if the discriminator is bounded away from zero and one and the derivatives of the generator and discriminator are bounded (this can be seen by straightforward differentiation). 
The last condition rules out degenerate cases; 
a prime example of such a degenerate case is the discriminator $D(x,\delta)$ being constant in $x$ for some $\delta$. 
Overall, these conditions are quite mild and would hold in many practical GAN applications.

Using the notation in Assumption GAN,
we can now restate the GAN minimax problem~(\ref{eq:IntroGANpopula})
from the Introduction as
\begin{equation}
\adjustlimits\inf_{\gamma\in\Gamma}\sup_{\delta\in\Delta}f(\gamma,\delta)
\qquad\text{with}\qquad 
f(\gamma,\delta)=\mathbb{E}[F(X,Z,\gamma,\delta)] .
\label{eq:GAN_population}
\end{equation}
In the traditional convex-concave setting ($f(\cdot,\delta)$ convex
for all fixed $\delta\in\Delta$ and $f(\gamma,\cdot)$ concave for
all fixed $\gamma\in\Gamma$), the classical von Neumann minimax theorem
implies that $\inf_{\gamma\in\Gamma}\sup_{\delta\in\Delta}f(\gamma,\delta)=\sup_{\delta\in\Delta}\inf_{\gamma\in\Gamma}f(\gamma,\delta)$
under mild conditions. In contrast to this, \citet{jin2020what} among
others have emphasized that in GAN applications $f$ is as a rule
non-convex and non-concave and the order in which minimization and
maximization are performed matters. Another point to note is that
in GAN applications the object of interest is not the optimal value
of problem (\ref{eq:GAN_population}), that is, 
\begin{equation}
V_{0}=\adjustlimits\inf_{\gamma\in\Gamma}\sup_{\delta\in\Delta}f(\gamma,\delta),\label{eq:V0}
\end{equation}
but rather the optimal solutions, the parameter values $(\gamma_{0},\delta_{0})$
that solve problem (\ref{eq:GAN_population}). Of main interest are
the $\gamma$-parameters appearing in the generator network as these
facilitate producing synthetic data according to $G(Z,\gamma_{0})$;
for completeness we consider also the $\delta$-parameters in the
discriminator network.

Let $\Theta_{0}\subseteq\Theta$ denote the set of optimal solutions
to (\ref{eq:GAN_population}). To describe $\Theta_{0}$, we introduce
the max-function 
\begin{equation}
\varphi(\gamma)=\sup_{\delta\in\Delta}f(\gamma,\delta).\label{eq:phi_func}
\end{equation}
A point $(\gamma_{0},\delta_{0})\in\Theta$ solves the GAN problem
(\ref{eq:GAN_population}) when it is a solution both to the inner
maximization problem (hence satisfying $f(\gamma_{0},\delta_{0})=
\sup_{\delta\in\Delta}f(\gamma_{0},\delta)=\varphi(\gamma_{0})$)
and to the outer minimization problem (hence satisfying
$\varphi(\gamma_{0})= \linebreak[3] \inf_{\gamma\in\Gamma}\varphi(\gamma)=V_{0}$).
Therefore the set of solutions $\Theta_{0}$ can be expressed as
\begin{equation}
\Theta_{0}=\{(\gamma_{0},\delta_{0})\in\Theta:f(\gamma_{0},\delta_{0})=\sup_{\delta\in\Delta}f(\gamma_{0},\delta)=\varphi(\gamma_{0})\text{ and }\varphi(\gamma_{0})=\inf_{\gamma\in\Gamma}\varphi(\gamma)=V_{0}\}.\label{eq:Theta0_A}
\end{equation}
Equivalently, $(\gamma_{0},\delta_{0})\in\Theta_{0}$ if and only
if $\max\{\varphi(\gamma_{0})-f(\gamma_{0},\delta_{0}),\varphi(\gamma_{0})-V_{0}\}=0$.
This motivates us to define the (population) criterion function
\begin{equation}
Q(\theta)=Q(\gamma,\delta)=\max\{\varphi(\gamma)-f(\gamma,\delta),\varphi(\gamma)-V_{0}\}.\label{eq:Q_func}
\end{equation}
The function $Q(\theta)$ is nonnegative for all $\theta$ and $\theta_{0}=(\gamma_{0},\delta_{0})\in\Theta_{0}$
if and only if $Q(\theta_{0})=0$. Therefore the set of solutions
(\ref{eq:Theta0_A}) can alternatively and concisely be characterized
as 
\begin{equation}
\Theta_{0}=\{\theta_{0}\in\Theta:Q(\theta_{0})=0\}.\label{eq:Theta0_B}
\end{equation}

Now consider the corresponding sample GAN minimax problem (\ref{eq:IntroGANsample}),
which can be written as 
\begin{equation}
\adjustlimits\inf_{\gamma\in\Gamma}\sup_{\delta\in\Delta}\hat{f}_{n}(\gamma,\delta)\qquad\text{with}\qquad\hat{f}_{n}(\gamma,\delta)=\frac{1}{n}\sum_{i=1}^{n}F(X_{i},Z_{i},\gamma,\delta)\label{eq:GAN_sample}
\end{equation}
and define the sample analogues of the quantities in (\ref{eq:V0})--(\ref{eq:Q_func}) as 
\begin{align}
\hat{V}_{n} 
	& =\adjustlimits\inf_{\gamma\in\Gamma}\sup_{\delta\in\Delta}\hat{f}_{n}(\gamma,\delta),\label{eq:VHat}\\
\hat{\varphi}_{n}(\gamma) 
	& =\sup_{\delta\in\Delta}\hat{f}_{n}(\gamma,\delta),\label{eq:phihat_func}\\
\hat{\Theta}_{n} 
	& =  \bigl\{  (\hat{\gamma}_{n},\hat{\delta}_{n})\in\Theta:    \; 
	\hat{f}_{n}(\hat{\gamma}_{n},\hat{\delta}_{n})  
	=  \sup_{\delta\in\Delta}\hat{f}_{n}(\hat{\gamma}_{n},\delta)
	=\hat{\varphi}_{n}(\hat{\gamma}_{n})   \text{ and } \label{eq:ThetaHat_A} \\
& \hspace*{90pt}  
	\hat{\varphi}_{n}(\hat{\gamma}_{n})
	=\inf_{\gamma\in\Gamma}\hat{\varphi}_{n}(\gamma)=\hat{V}_{n}  \bigr  \},\notag \\
\hat{Q}_{n}(\theta) 
	& =\hat{Q}_{n}(\gamma,\delta)
	=\max\{\hat{\varphi}_{n}(\gamma)-\hat{f}_{n}(\gamma,\delta),
	\hat{\varphi}_{n}(\gamma)-\hat{V}_{n}\}.\label{eq:Qhat_func}
\end{align}
These quantities have interpretations similar to their population
counterparts: $\hat{V}_{n}$ is the optimal value of the sample GAN
problem (\ref{eq:GAN_sample}), $\hat{\varphi}_{n}(\gamma)$ is a
sample max-function, $\hat{\Theta}_{n}$ denotes the set of (exact)
solutions to the sample GAN problem (\ref{eq:GAN_sample}), and $\hat{Q}_{n}(\theta)$
is a (nonnegative) sample criterion function that allows us to express
the set $\hat{\Theta}_{n}$ as 
\begin{equation}
\hat{\Theta}_{n}=\{\hat{\theta}_{n}\in\Theta:\hat{Q}_{n}(\hat{\theta}_{n})=0\}.\label{eq:ThetaHat_B}
\end{equation}

Finding a solution to the sample GAN problem (\ref{eq:GAN_sample}),
let alone the entire set of optimal solutions $\hat{\Theta}_{n}$,
is a challenging task. From a practical machine learning perspective,
this `training of the GAN' is of principal interest and extensive
research efforts have been made to devise algorithms to do this (see
the references listed in the Introduction). In this paper, we do not
consider these algorithms but rather just assume that the sample GAN
problem is solved using some available method. These algorithms typically
search for approximate rather than exact solutions:\ 
for some small non-negative constant $\tau$, approximate solutions to the
sample GAN problem (\ref{eq:GAN_sample}) are points $(\hat{\gamma}_{n},\hat{\delta}_{n})\in\Gamma\times\Delta$
satisfying
\[
\hat{f}_{n}(\hat{\gamma}_{n},\hat{\delta}_{n})\geq\sup_{\delta\in\Delta}\hat{f}_{n}(\hat{\gamma}_{n},\delta)-\tau\quad\text{and}\quad\hat{\varphi}_{n}(\hat{\gamma}_{n})\leq\inf_{\gamma\in\Gamma}\hat{\varphi}_{n}(\gamma)+\tau.
\]
Such points $(\hat{\gamma}_{n},\hat{\delta}_{n})$ approximately solve
both the inner maximization problem and the outer minimization problem,
with $\tau$ determining the slackness allowed in these maximization
and minimization problems. Somewhat more generally, let $\tau_{n}$
be a sequence of non-negative random variables such that $\tau_{n}\overset{p}{\to}0$
(where $\overset{p}{\to}$ denotes convergence in probability).
Define in analogy with the above 
\begin{align}
\hat{\Theta}_{n}(\tau_{n})
& =\bigl \{  (\hat{\gamma}_{n},\hat{\delta}_{n})\in\Theta:    \;
	\hat{f}_{n}(\hat{\gamma}_{n},\hat{\delta}_{n})
	\geq\sup_{\delta\in\Delta}\hat{f}_{n}(\hat{\gamma}_{n},\delta)-\tau_{n}  \text{ and }  \label{eq:ThetaHatTau_A}  \\
& \hspace*{90pt}  
	\hat{\varphi}_{n}(\hat{\gamma}_{n})
	\leq\inf_{\gamma\in\Gamma}\hat{\varphi}_{n}(\gamma)+\tau_{n}  \bigr  \}  \notag
\end{align}
as the set of approximate solutions to the sample GAN problem (\ref{eq:GAN_sample}).
Noting that the two inequalities in (\ref{eq:ThetaHatTau_A}) can
be expressed as $\hat{\varphi}_{n}(\hat{\gamma}_{n})-\hat{f}_{n}(\hat{\gamma}_{n},\hat{\delta}_{n})\leq\tau_{n}$
and $\hat{\varphi}_{n}(\hat{\gamma}_{n})-\hat{V}_{n}\leq\tau_{n}$,
respectively, and recalling the definition of the sample criterion
function $\hat{Q}_{n}(\theta)$ in (\ref{eq:Qhat_func}) allows us
to characterize the set $\hat{\Theta}_{n}(\tau_{n})$ as
\begin{equation}
\hat{\Theta}_{n}(\tau_{n})=\{\hat{\theta}_{n}\in\Theta:\hat{Q}_{n}(\hat{\theta}_{n})\leq\tau_{n}\}.\label{eq:ThetaHatTau_B}
\end{equation}
Setting $\tau_{n}=0$ one obtains as a special case the set of exact
solutions, $\hat{\Theta}_{n}=\hat{\Theta}_{n}(0)$. 
In the next section we consider 
the properties of $\hat{\Theta}_{n}(\tau_{n})$ as an estimator of $\Theta_{0}$.

Before proceeding, a remark about the interpretation of $\Theta_{0}$ is in place.
In GAN applications, it is typically not realistic to assume that the generated synthetic data would perfectly mimic the real observations. 
In line with this, the GAN formulation allows for misspecification:
It is not assumed that $G(Z,\gamma)$ would for some $\gamma$ have the same distribution as $X$ does.
In this sense, the elements of $\Theta_{0}$ do not correspond to any `true' parameter values. 
Nevertheless, an interpretation can be given. 
To do so, momentarily consider a somewhat idealized version of the GAN problem where in (1) 
the supremum over $\delta$ and the parametric class of discriminator functions $D(\cdot,\delta)$ (from $\mathcal{X}$ to $(0,1)$)
is replaced with the supremum over all measurable functions $D(\cdot)$ from $\mathcal{X}$ to $(0,1)$. 
It turns out that in this case the GAN problem in (1) reduces to the minimization problem 
$\inf_{\gamma\in\Gamma} 2 [ \text{JSD}(P_X,P_{\gamma}) -  \ln 2 ]$,
where $P_X$ denotes the probability distribution of $X$, 
$P_{\gamma}$ the probability distribution of $G(Z,\gamma)$, 
and $\text{JSD}(P_X,P_{\gamma})$ the so-called Jensen-Shannon divergence between these two  distributions
(see \citealp[Sec 2]{biau2020some}, and \citealp[Secs 1--2]{belomestny2021rates},
for the technical details and the definition of the Jensen-Shannon divergence).
This yields an idealized interpretation of the GAN problem:
heuristically, one can interpret the elements of $\Theta_{0}$ to correspond to 
$\gamma_0$ that minimize the Jensen-Shannon divergence 
between the true target distribution $P_X$ 
and the generated distribution of $G(Z,\gamma_0)$
(of course, this interpretation is not entirely accurate as in practice
the discriminators employed are parametric neural networks). 
Note also that this interpretation is akin to conventional maximum likelihood (ML) estimation in misspecified models,
where (under appropriate assumptions) the ML estimator converges to a parameter, say again $\theta_0$, 
that minimizes the Kullback-Leibler divergence between the true distribution and the distribution corresponding to parameter $\theta_0$
(see, e.g., \citealp[Ex 5.25]{vanderVaart1998asymptotic}).

\section{Consistent estimation}

In GAN minimax problems both $\hat{\Theta}_{n}(\tau_{n})$ and $\Theta_{0}$
are typically set-valued and not singletons and an appropriate notion
of distance between sets is required. One commonly used generalization
of the Euclidean distance $|\cdot|$ is the Hausdorff distance (see,
e.g., \citealp[Sec 4.C]{rockafellar2009variational}, or \citealp[Appendix D]{molchanov2017theory}).
For any two non-empty bounded subsets $A$ and $B$ of some Euclidean
space, the Hausdorff distance between $A$ and $B$ is defined as
\[
d_{H}(A,B)=\max\bigl\{\sup_{a\in A}d(a,B),\sup_{b\in B}d(b,A)\bigr\},
\]
where $d(a,B)=\inf_{b\in B}|a-b|$ is the shortest distance from the
point $a$ to the set $B$. That is, the Hausdorff distance is the
greatest distance from an arbitrary point in one of the sets to the
closest neighboring point in the other set. The Hausdorff distance
$d_{H}$ is a metric for the family of non-empty compact sets and
for such sets $d_{H}(A,B)=0$ if and only if $A=B$.

The consistency result we aim to prove is Hausdorff consistency in
the sense that\linebreak[4] 
$d_{H}(\hat{\Theta}_{n}(\tau_{n}),\Theta_{0})\overset{p}{\to}0$. Establishing
this requires us to show that both of the `one-sided Hausdorff consistency'
conditions 
\begin{equation}
\text{(a)}\;  \sup_{\theta\in\hat{\Theta}_{n}(\tau_{n})}d(\theta,\Theta_{0})\overset{p}{\to}0
\qquad\text{and}\qquad
\text{(b)}\;  \sup_{\theta\in\Theta_{0}}d(\theta,\hat{\Theta}_{n}(\tau_{n}))\overset{p}{\to}0\label{eq:1SidedConsistency}
\end{equation}
hold. Heuristically, the former condition in (\ref{eq:1SidedConsistency}) guarantees that $\hat{\Theta}_{n}(\tau_{n})$
is not too large compared to $\Theta_{0}$, whereas the latter condition
ensures $\hat{\Theta}_{n}(\tau_{n})$ is large enough to cover all
of $\Theta_{0}$. Establishing (\ref{eq:1SidedConsistency}a)
follows the pattern of a standard consistency proof and relies on
a suitable uniform law of large numbers combined with an appropriate
set-identification condition for $\Theta_{0}$. Proving (\ref{eq:1SidedConsistency}b) relies on also knowing
the rate of this uniform convergence. The following assumption formally
states the required high-level conditions; here $O_{p}(1)$ stands
for a sequence of random variables that is bounded in probability.
\begin{assumption}
\mbox{}
\begin{lyxlist}{(a)}
\item [{(a)}] $\sup_{\theta\in\Theta}|\hat{f}_{n}(\theta)-f(\theta)|\overset{p}{\to}0$
with the function $f(\theta)$ continuous in $\theta$.
\item [{(b)}] Part (a) holds with $\sup_{\theta\in\Theta}n^{1/2}|\hat{f}_{n}(\theta)-f(\theta)|=O_{p}(1)$.
\end{lyxlist}
\end{assumption}
The high-level conditions in this assumption can be verified using
various sets of sufficient conditions. For instance, Assumption
1(a) holds when Assumption GAN is combined with the mild moment condition
$\mathbb{E}[\sup_{\theta\in\Theta}|F(X,Z,\theta)|]<\infty$, and 1(b)
holds under Assumption Smooth GAN (justifications for these statemens
are given in Appendix C). Assumption 1 thus holds in most practical
GAN applications.

Our consistency results also require certain conditions for the slackness
sequence $\tau_{n}$.
\begin{assumption}
\mbox{}
\begin{lyxlist}{(a)}
\item [{(a)}] $\tau_{n}$ is a sequence of non-negative random variables
such that $\tau_{n}\overset{p}{\to}0$.
\item [{(b)}] $\tau_{n}$ is a sequence of positive random variables
such that $\tau_{n}\overset{p}{\to}0$ and \linebreak[3] $n^{-1/2}/\tau_{n}\overset{p}{\to}0$.
\end{lyxlist}
\end{assumption}
Part (a) of Assumption 2 allows for the possibility that $\tau_{n}$
is identically zero while this is ruled out in part (b). In (b) it
is additionally assumed that the convergence of $\tau_{n}$ to zero
is slower than that of $n^{-1/2}$ in the sense that $n^{-1/2}/\tau_{n}\overset{p}{\to}0$;
for instance, $\tau_{n}=n^{-0.49}$ is a possibility. Of course, (b)
implies (a).

We can now state our consistency theorem for the solutions of the
GAN minimax problem (the proof is given in Appendix A). 
\begin{thm}
Suppose Assumption GAN holds. 
\begin{lyxlist}{(a)}
\item [{(a)}] Assumptions 1(a) and 2(a) imply that (\ref{eq:1SidedConsistency}a) holds
(and also that $\hat{V}_{n}\overset{p}{\to}V_{0}$).
\item [{(b)}] Assumptions 1(b) and 2(b) imply that (\ref{eq:1SidedConsistency}b) holds 
so that $d_{H}(\hat{\Theta}_{n}(\tau_{n}),\Theta_{0})\overset{p}{\to}0$.
\end{lyxlist}
\end{thm}
Part (a) of Theorem 1 establishes the former one-sided Hausdorff consistency
condition in (\ref{eq:1SidedConsistency}); for completeness, consistency
of the optimal value $\hat{V}_{n}$ is also given. Under the stronger
conditions in part (b), the latter condition in (\ref{eq:1SidedConsistency})
is also obtained, thus establishing the desired Hausdorff consistency
result. 

Theorem 1 is based on the Hausdorff consistency results of \citet[Thm 3.1]{chernozhukov2007estimation}
in general set-identified situations. \citet[Thm 5.9]{shapiro2009lectures}
give (an almost sure version) of part (a) of Theorem 1 in a general
minimax setting. In the GAN setting, a consistency result in the case
of a unique solution is given in \citet[Thm 4.2]{biau2020some}. (Note
that when $\Theta_{0}$ is a singleton-set consisting of one point
$\theta_{0}$, Theorem 1(a) shows that any element of $\hat{\Theta}_{n}(\tau_{n})$
is consistent for $\theta_{0}$.) We are not aware of consistency
results in the GAN setting in the case of multiple solutions and the
result of Theorem 1 is reassuring in that estimation in the GAN setting
will be consistent regardless of the (potentially infinite) number
of solutions. The one-sided consistency result of Theorem 1(a) covers
the case of exact solutions ($\tau_{n}=0$), whereas the two-sided
case in part (b) requires us to (somewhat arbitrarily) choose a strictly
positive slackness sequence $\tau_{n}$. Such a choice is not needed
in the procedure for forming confidence sets for $\Theta_{0}$ that
we consider next.

\section{Confidence sets}

A confidence set for 
the GAN minimax solutions
$\Theta_{0}$ is a random set that covers the
entire $\Theta_{0}$ with a prespecified probability. 
(As in GANs the individual elements $\theta_{0}$ of $\Theta_{0}$ 
do not have any special interpretation attached to them, 
a confidence set covering all of $\Theta_{0}$ 
and not just some particular element of it is appropriate.)
Let $1-\alpha$ denote the desired coverage probability
(such as $95\%$) where $\alpha\in(0,1)$. We aim to construct confidence
sets $\hat{\CS}_{n,1-\alpha}$ that contain $\Theta_{0}$ 
with at least probability $1-\alpha$,
\begin{equation}
\liminf_{n\to\infty}\mathbb{P}[\Theta_{0}\subseteq\hat{\CS}_{n,1-\alpha}]\geq1-\alpha,\label{eq:CS_level_1-a}
\end{equation}
that is, sets that are conservatively asymptotically consistent at
level $1-\alpha$. (Here $\mathbb{P}$ refers to the probability measure
of $(X,Z)$ which we consider fixed.) 

In the traditional point-identified setting in which $\Theta_{0}$
consists of a single point $\theta_{0}$, a conventional route to
forming confidence sets is to consider a Taylor approximation of a
certain function around $\theta_{0}$; this leads to the (often
Gaussian) distribution of an appropriately rescaled estimator from
which confidence sets for $\theta_{0}$ are then obtained in a straightforward
manner. \citet{biau2020some} consider this point-identified case
and, focusing on the generator network $\gamma$-parameters, in their
Theorem 4.3 state that $n^{1/2}(\hat{\gamma}_{n}-\gamma_{0})$ converges
in distribution to a Gaussian random variable. When $\Theta_{0}$
is not a singleton set alternative approaches are called for. 
Given the results of \citet{biau2020some}, we
focus only on the case of multiple (more than one) solutions. This
is the typical case in GAN applications. Our approach follows the
recent developments in the partial identification literature in econometrics;
we rely in particular on the results in \citet{chernozhukov2007estimation}
and \citet{romano2010inference} 
(for further related references, see the survey of this literature in \citealp{molinari2020microeconometrics}).

The confidence sets we consider are based on appropriate lower contour
sets of the rescaled criterion function $n^{1/2}\hat{Q}_{n}(\theta)$.
To motivate the subsequent technical developments, we begin with some informal remarks. 
First, if the $1-\alpha$ quantile $c_{n,1-\alpha}$ of the distribution of
$\sup_{\theta\in\Theta_{0}}n^{1/2}\hat{Q}_{n}(\theta)$ was known,
one could form the (infeasible) confidence set 
$\CS_{n,1-\alpha} = \{\theta\in\Theta:n^{1/2}\hat{Q}_{n}(\theta)\leq c_{n,1-\alpha}\}$
that satisfies $\mathbb{P}[\Theta_{0}\subseteq\CS_{n,1-\alpha}]\geq1-\alpha$
(because $\Theta_{0}\subseteq\CS_{n,1-\alpha}$ if and only if $\sup_{\theta\in\Theta_{0}}n^{1/2}\hat{Q}_{n}(\theta)\leq c_{n,1-\alpha}$).
Second, consider the statistic $\sup_{\theta\in S}n^{1/2}\hat{Q}_{n}(\theta)$
where $S$ is some nonempty subset of $\Theta$; 
under appropriate conditions, it can be shown that as $n\to\infty$,
$\sup_{\theta\in S}n^{1/2}\hat{Q}_{n}(\theta)$ remains stochastically
bounded for $S\subseteq\Theta_{0}$ and diverges to infinity for $S\nsubseteq\Theta_{0}$. 
Now, to form our confidence sets
we approximate the $1-\alpha$ quantile of $\sup_{\theta\in S}n^{1/2}\hat{Q}_{n}(\theta)$
for various sets $S$ using a suitable resampling method, 
and then locate a confidence set for $\Theta_{0}$ 
based on the different behavior of $\sup_{\theta\in S}n^{1/2}\hat{Q}_{n}(\theta)$
in the two cases $S\subseteq\Theta_{0}$ and $S\nsubseteq\Theta_{0}$. 
As our situation involves non-standard features, resampling based on 
standard bootstrap is not appropriate and instead we use a procedure 
based on subsampling (see \citealp{politis1999subsampling}).

To formalize this discussion, we next introduce some notation.
For any nonempty subset $S$ of $\Theta$, let $L_{n}(S,x)$ denote the cumulative
distribution function (cdf) of the statistic $\sup_{\theta\in S}n^{1/2}\hat{Q}_{n}(\theta)$
and $c_{n,1-\alpha}(S)$ the corresponding (smallest) $1-\alpha$
quantile, that is, 
\begin{align*}
L_{n}(S,x)&=\mathbb{P}[\sup_{\theta\in S}n^{1/2}\hat{Q}_{n}(\theta)\leq x]\quad\text{and} \\
c_{n,1-\alpha}(S)&=\inf\{x\in\mathbb{R}:L_{n}(S,x)\geq1-\alpha\}.
\end{align*}
Let $b$ denote the subsample size that is assumed to satisfy the
usual requirements $b\to\infty$ and $b/n\to0$ as $n\to\infty$;
for notational simplicity, the dependence of $b$ on $n$ is not made
explicit. The number of different subsamples is denoted by $N_{n}$
($=n!/(b!(n-b)!)$) and the subsample statistics $\hat{Q}_{n,b,i}(\theta)$
($i=1,\ldots,N_{n}$) are defined exactly as $\hat{Q}_{n}(\theta)$
in (\ref{eq:Qhat_func}) but based on the $i$th subsample
of size $b$ rather than the full sample. The subsampling counterparts
of $L_{n}(S,x)$ and $c_{n,1-\alpha}(S)$ are $\hat{L}_{n,b}(S,x)$,
the empirical distribution function of the 
(centered) subsample statistics 
$\{ b^{1/2}[ \sup_{\theta\in S}\hat{Q}_{n,b,i}(\theta) - \sup_{\theta\in S}\hat{Q}_{n}(\theta) ]
:i=1,\ldots,N_{n}\}$,
and $\hat{c}_{n,b,1-\alpha}(S)$, the corresponding $1-\alpha$ sample
quantile, defined as 
\begin{align}
\hat{L}_{n,b}(S,x) & = 
\frac{1}{N_{n}}\sum_{i=1}^{N_{n}}
1\Bigl(  
b^{1/2}[ \sup_{\theta\in S}\hat{Q}_{n,b,i}(\theta) - \sup_{\theta\in S}\hat{Q}_{n}(\theta) ]
  \leq x  \Bigr)
\quad\text{and} \label{eq:SS_ecdf}  \\
\hat{c}_{n,b,1-\alpha}(S) & =\inf\{x\in\mathbb{R}:\hat{L}_{n,b}(S,x)\geq1-\alpha\}. \notag
\end{align}

We can now state the iterative procedure that we use to construct
the desired confidence sets. This procedure is akin to step-down procedures
used in multiple hypothesis testing problems and follows 
\citet[Sec 9.1]{lehmann2005testing} and \citet{romano2010inference}. 

\begin{myprocedure}
Set $S_{1}=\Theta$ and for $j=1,2,\ldots$ do the following: If 
\[ \sup_{\theta\in S_j}n^{1/2}\hat{Q}_{n}(\theta)\leq\hat{c}_{n,b,1-\alpha}(S_{j}), \] 
then set $\hat{\CS}_{n,1-\alpha}=S_{j}$ and stop; 
otherwise, set 
$S_{j+1}=\{\theta\in\Theta:n^{1/2}\hat{Q}_{n}(\theta)\leq\hat{c}_{n,b,1-\alpha}(S_{j})\}$ 
and continue iteration.
\end{myprocedure}

This procedure starts from the full parameter space $S_{1}=\Theta$
and iteratively discards $\theta$'s from the $S_{j}$-sets until
a suitable confidence set is formed. Overall, the procedure is certainly
computer-intensive yet feasible (to reduce computational costs, one
can use just a subset of the $N_{n}$ subsamples without affecting
the validity of our results; cf.~\citealp[Cor 2.4.1]{politis1999subsampling}).
Implementation of the procedure is illustrated 
in a small numerical example in the next section.
Theorem 2 below shows that, under appropriate assumptions, the confidence
sets $\hat{\CS}_{n,1-\alpha}$ formed using Procedure 1 are valid
confidence sets in the sense that (\ref{eq:CS_level_1-a}) holds.
The key requirement for this subsampling-based procedure to work is
the following high-level assumption.

\begin{assumption}
The statistic $\sup_{\theta\in\Theta_{0}}n^{1/2}\hat{Q}_{n}(\theta)$
converges in distribution to some limiting random variable $L$.
\end{assumption}

The convergence requirement in Assumption 3 is a standard high-level
condition needed for the validity of subsampling procedures 
(see \citealp[Sec 2.2]{politis1999subsampling}).
Ensuring that Assumption 3 holds in 
GAN minimax problems 
turns out to be
particularly challenging and we resort to empirical process theory
for this. (A thorough account of the needed empirical process theory
can be found in the monograph of \citet{vandervaart1996weak}.)
With the set $\Theta$ compact as before, let $l^{\infty}(\Theta)$ denote
the space of all bounded real-valued functions on $\Theta$ equipped
with the supremum norm and $C(\Theta)$ stand for the subspace of
those functions that are continuous. Random functions such as $\hat{f}_{n}(\cdot)$
 are viewed as maps from appropriate underlying
probability spaces to $l^{\infty}(\Theta)$, and $\rightsquigarrow$
denotes weak convergence as defined in \citet[Sec 1.3]{vandervaart1996weak}.

A convenient and broadly applicable high-level assumption that we
make is that a suitable functional central limit theorem (Donsker
property) holds.

\begin{assumption}
Suppose $n^{1/2}(\hat{f}_{n}-f)\rightsquigarrow\mathbb{G}$ in $l^{\infty}(\Theta)$
with $f\in C(\Theta)$ and $\mathbb{G}$ a tight mean-zero Gaussian
process taking values in $C(\Theta)$ with probability one and such
that $\inf_{\theta\in\Theta_{0}}\mathbb{E}[(\mathbb{G}(\theta))^{2}]>0$.
\end{assumption}

Assumption 4 is very general and rather technical. 
Importantly, we note that Assumption Smooth GAN implies the
validity of Assumption 4 (justification of this is given in Appendix
C). Alternatively, Assumption 4 also holds under much weaker conditions and can
be verified using a variety of different methods in empirical process
theory (for details, see \citealp{vandervaart1996weak}).

Showing that Assumption 3 follows from Assumption 4 (and other additional
conditions) requires rather long and technical details.
In order to not distract from the main issue, we outline the key points
of the argument here and relegate the technical details to Appendix
B. First, it is shown that for a suitably defined map $\phi:l^{\infty}(\Theta)\to\mathbb{R}$,
the statistic $\sup_{\theta\in\Theta_{0}}n^{1/2}\hat{Q}_{n}(\theta)$
in Assumption 3 can be expressed as $\sup_{\theta\in\Theta_{0}}n^{1/2}\hat{Q}_{n}(\theta)=n^{1/2}(\phi(\hat{f}_{n})-\phi(f))$.
Second, the map $\phi$ is shown to be (directionally) differentiable
in an appropriate sense
-- this step is the technically most delicate one (and Lemma 3 in Appendix B containing the 
key novel results may be of independent interest). 
Third, the previous facts enable us to apply
a particular version of the functional delta method to deduce that
$n^{1/2}(\phi(\hat{f}_{n})-\phi(f))\rightsquigarrow\phi_{f}^{\prime}(\mathbb{G})$
where $\phi_{f}^{\prime}$ denotes a certain (directional) derivative
of the mapping $\phi$ at $f$ (and $\mathbb{G}$ is the Gaussian
process in Assumption 4). Thus Assumption 3 follows.
To prove the differentiability result mentioned,
additional assumptions are needed and we focus on the following leading
case. 
\begin{assumption}
The set $\Gamma_{0}=\{\gamma_{0}\in\Gamma:(\gamma_{0},\delta_{0})\in\Theta_{0}\text{ for some }\delta_{0}\in\Delta\}$
is finite.
\end{assumption}
Assumption 5 requires that the outer minimization problem in the population
minimax problem (\ref{eq:GAN_population}) is solved at a finite number
of $\gamma$-values only (the inner maximization problem is allowed
to have infinitely many solutions). This assumption is quite reasonable
in practical GAN applications and is needed to prove the mentioned
differentiability result. 

All the details of the preceding discussion are available in Appendix B, 
where we prove the following lemma. 
\begin{lem}
Suppose Assumptions GAN, 4, and 5 hold.
Then Assumption 3 holds with the limiting distribution
of\/ $\sup_{\theta\in\Theta_{0}}n^{1/2}\hat{Q}_{n}(\theta)$ given
in expression (\ref{eq:LimitingRV}) in Appendix B. 
\end{lem}
Lemma 1 offers one convenient way to verify that Assumption 3 holds
in the GAN setting. Even with the finiteness requirement of Assumption
5, the technical details in Appendix B are rather long and the limiting
distribution in (\ref{eq:LimitingRV}) quite complicated. 
Assumption 3 could potentially be verified using weaker conditions 
but we do not pursue this.

After these preparations, we are now ready to state our 
main result regarding the confidence sets formed by Procedure 1. 
Let $L(\cdot)$ denote the cdf of the limiting random variable $L$ in Assumption 3,
$c_{1-\alpha}=\inf\{x\in\mathbb{R}:L(x)\geq1-\alpha\}$ 
the corresponding $1-\alpha$ quantile, and 
$\lim_{c\smalluparrow c_{1-\alpha}}L(c)$ 
the limit from the left of $L(\cdot)$ at $c_{1-\alpha}$.
The following theorem is proved in Appendix B.

\begin{thm}
Suppose Assumptions GAN and 3 hold, 
the population GAN minimax problem (\ref{eq:GAN_population}) has multiple solutions, 
the confidence sets $\hat{\CS}_{n,1-\alpha}$ ($n=1,2,\ldots$)  are formed using Procedure 1, and  the subsampling size $b$ is such that 
$b\to\infty$ and $b/n\to0$ as $n\to\infty$. 
Then the confidence sets $\hat{\CS}_{n,1-\alpha}$ satisfy 
\[
\liminf_{n\to\infty}\mathbb{P}[\Theta_{0}\subseteq\hat{\CS}_{n,1-\alpha}]\geq\lim_{c\smalluparrow c_{1-\alpha}}L(c).
\]
\end{thm}

Theorem 2 shows that the confidence sets $\hat{\CS}_{n,1-\alpha}$
have an asymptotic coverage probability of at least $1-\alpha$
when the cdf $L(\cdot)$ is continuous at $c_{1-\alpha}$
(as then $\lim_{c\smalluparrow c_{1-\alpha}}L(c)=1-\alpha$; 
if $L(\cdot)$ has a jump at $c_{1-\alpha}$, then 
$\lim_{c\smalluparrow c_{1-\alpha}}L(c)<1-\alpha\leq L(c_{1-\alpha})$).
In general the limiting distribution $L$ is rather complicated and
verifying the continuity of $L(\cdot)$ requires some more concrete
assumptions. In Appendix B we show that this continuity holds for instance
when $\Theta_{0}$ is a finite set of the form $\Theta_{0}=\{(\gamma_{01},\delta_{01}),\ldots,(\gamma_{0K},\delta_{0K})\}$
for some finite $K>1$ and with distinct $\gamma$'s and $\delta$'s
(i.e., $\gamma_{0i}\neq\gamma_{0j}$ and $\delta_{0i}\neq\delta_{0j}$
for all $i\neq j$). 
This case with multiple but finite number of
distinct solutions covers many practical GAN applications
-- in particular, in this case Assumption Smooth GAN
suffices for the validity of Theorem 2.

In previous work, confidence regions for general set-identified parameters
were considered by \citet{chernozhukov2007estimation} and \citet{romano2010inference}.
Our Theorem 2 is based on their results, especially on Theorem 2.2
of \citet{romano2010inference}. 
However, these previous works do not consider minimax problems.
\citet[Thm 3.1]{shapiro2008asymptotics}
and \citet[Thm 5.10]{shapiro2009lectures} study the asymptotic distribution
of the optimal value in convex-concave minimax problems but do not
consider inference for the solutions of minimax problems. 
The confidence sets we provide
are novel in the context of general minimax problems.
\citet[Thm 4.3]{biau2020some}
consider the GAN problem assuming a single unique solution exists
and give a result on the asymptotic normality of the generator network
parameters; however, they remark (pp.~1560--1561) that the assumed
uniqueness is ``hardly satisfied in the high-dimensional context
of (deep) neural networks'' and call for generalizations.
Theorem 2 builds on these previous results and provides the first inference
procedure for GAN parameters
in the presence of multiple solutions.
Our construction of the confidence sets is specifically tailored for the case of multiple solutions; in the case of a single solution, 
confidence sets can be formed based on the results of \citet{biau2020some} 
(for details, see Appendix B).

The confidence sets of this section
facilitate statistical inference for the solutions of  
GAN minimax problems.
Regarding the potential practical use of these confidence sets
in GAN applications, recall from the Introduction that when
$(\gamma_{0},\delta_{0})\in\Theta_{0}$, the generator $G(Z,\gamma_{0})$
gives the researcher a mechanism to produce synthetic data. In recent
work, \citet[Sec 3.2, Fig 4]{karras2020style} have considered the
effect of small variations in the input noise on the synthetic images
produced. Similarly, one could consider how the produced images (or,
more generally, the distribution of $G(Z,\gamma_{0})$) change when
the generator parameters 
vary, say, within a certain confidence set.
Related to this, the effect of varying the generator parameters
(through parameter interpolation) 
has recently been considered in 
\citet{wang2019deep} and \citet[Secs 4.2 and 5]{pan2022exploiting}.
Exploring the effect of such variations on the synthetic data produced 
in other GAN applications 
(such as DNA sequences, art, or music) would also be interesting.
We leave these issues for future research.

\section{Numerical illustration}

We next illustrate our results in small numerical experiments.
We consider a toy example that is as simple as possible, thus allowing us to both
analytically solve the population GAN problem as well as 
to graphically illustrate consistent estimation and confidence sets for $\Theta_0$ in two dimensions.
It should be emphasized that real GAN applications used in practice are of course remarkably more complicated than our toy illustration.

In our example, we use a GAN to mimic data from a more complicated univariate distribution by simple Gaussian noise. The set-up is quite similar to one of the examples \citet[p.\ 1552--]{biau2020some} use to illustrate their results.
The real data $X$ is assumed to follow the so-called claw distribution (\citealp{marron1992exact}): 
Letting $N(x;\mu,\sigma^2)$ denote the probability density function (pdf) of a normal random variable with mean $\mu$ and variance $\sigma^2$, the claw density is defined as 
\[ p_{\text{claw}}(x) = \tfrac{1}{2} N(x;0,1)  + \textstyle{\displaystyle{\sum}_{k=0}^4} \tfrac{1}{10} N(x;\tfrac{k}{2} -1,0.01);
 \]
this pdf together with the standard normal one are illustrated in the top left graph of Figure 1.
The random noise $Z$ is assumed to be standard normally distributed, and the generator
function is a simple shift formulated as 
\[
G(z, \gamma ) = z + \cos(\gamma \pi)    
\quad\text{with}\quad \gamma\in[0,2].
\]
With this very simple formulation, the distribution of $G(Z, \gamma )$ never matches that of $X$ but is close to it when $\gamma = 0.5$ or $\gamma = 1.5$ (corresponding to no shift at all). 
The cosine function is used to incorporate multiple solutions in a transparent manner.
The discriminator we employ is a simple one parameter function
(motivated by \citealp[eqn (2)]{goodfellow2014generative}, and 
\citealp[eqn (2.1)]{biau2020some})
given by
\[
D(x, \delta ) = \frac{p_{\text{claw}}(x)}{p_{\text{claw}}(x) + N(x;\cos(\delta \pi),1)}
\quad\text{with}\quad \delta\in[0,2].
\]
Intuitively, the two densities $p_{\text{claw}}(x)$ and $N(x;\cos(\delta \pi),1)$ are close when $\delta = 0.5$ or $\delta = 1.5$.

More formally, in this toy example the population GAN problem (\ref{eq:GAN_population}) can be easily solved (see \citealp[Sec 2]{biau2020some}, and \citealp[Secs 1--2]{belomestny2021rates},
for further helpful details):
For any fixed $\gamma$, the inner maximization problem in (\ref{eq:GAN_population}) is solved for those $\delta$ that satisfy $\cos(\delta \pi) = \cos(\gamma \pi)$ as in these cases the discriminator $D(x, \delta )$ coincides with the so-called optimal discriminator; therefore problem (\ref{eq:GAN_population}) reduces to minimizing the Jensen-Shannon divergence between the probability distributions of $X$ and $G(Z, \gamma )$ (cf. the last paragraph of Section 2) which happens for $\gamma = 0.5$ or $\gamma = 1.5$. Therefore the population GAN problem has four solutions and
$\Theta_0 = \{ (0.5,0.5) , (0.5,1.5) , (1.5,0.5) , (1.5,1.5) \}$.
These are the four points that solve 
$\inf_{\gamma\in\Gamma}\sup_{\delta\in\Delta}f(\gamma,\delta)$ in (\ref{eq:GAN_population})
or, alternatively, that solve $Q(\gamma,\delta)=0$ (see (\ref{eq:Q_func}) and (\ref{eq:Theta0_B})).
The second and third plots in the first row of Figure 1 illustrate 
functions $f(\gamma,\delta)$ (second plot) and $Q(\gamma,\delta)$ (third plot) 
as well as the four solutions in $\Theta_0$ (the four red dots in these plots).
Note that already in this toy example, the landscapes of $f$ and $Q$ are somewhat non-trivial.

Now consider the sample GAN problem (\ref{eq:GAN_sample}) and consistency of the estimator
$\hat{\Theta}_{n}(\tau_n)$. 
Taking $n$\ IID draws $X_1,\ldots,X_n$ and $Z_1,\ldots,Z_n$ from $X$ and $Z$ as above, the
second row of Figure 1 illustrates $\hat{\Theta}_{n}(\tau_n)$ for four different choices of $n$ and $\tau_n$: 
$n=100$ with $\tau_n=0.1n^{-0.25}$ and  $\tau_n=0.1n^{-0.49}$ (first two plots) and 
$n=10.000$ with $\tau_n=0.1n^{-0.49}$ and $\tau_n=0$ (last two plots). 
Here and in all other calculations, we consider values of
$\gamma$ and $\delta$ over the grid $\{0,0.01, \ldots ,1.99,2\}$.
Comparing the first two estimators demonstrates the effect of a slower vs faster convergence of $\tau_n$ to zero, while the effect of an increasing sample size is seen by comparing the second estimator with the third. The fourth plot illustrates that the exact estimator $\hat{\Theta}_{n}(0)$ is not necessarily Hausdorff consistent as it may not contain all four elements of $\Theta_0$. 

We next illustrate forming confidence sets using Procedure 1: 
One starts from the full parameter space $S_1 = \Theta$ and 
iteratively discards $\theta$'s from the sets $S_j$ 
based on the appropriate quantile of suitable subsample statistics. 
The third and fourth rows of Figure 1 illustrate the working of Procedure 1 in one 
simulated data set, where we used sample size $n=1000$, 
subsample size $b=501$, $\alpha=0.20$, 
and $200$ randomly chosen subsamples
(as mentioned after Procedure 1, one can use only a subset of the possible subsamples and not all of them; note that the total number of possible subsamples in our exercise is astronomical as $\binom{1000}{501}\approx 10^{299}$).
The third row of Figure 1 plots the sets $S_2$, $S_3$, $S_4$, and 
the final confidence set $\hat{\CS}_{n,1-\alpha}=S_6$.
Below these are shown (for $j=2,3,4,6$) the empirical distribution function $\hat{L}_{n,b}(S_j,x)$ together with the  $1-\alpha$ sample quantile $\hat{c}_{n,b,1-\alpha}(S_{j})$ (red dotted line) and the quantity $\sup_{\theta\in S_j}n^{1/2}\hat{Q}_{n}(\theta)$ (blue dotted line).
In this particular data set, on round 6 of Procedure 1, 
$\sup_{\theta\in S_6}n^{1/2}\hat{Q}_{n}(\theta) < \hat{c}_{n,b,1-\alpha}(S_6)$ 
and thus the procedure stops and the confidence set is formed as 
$\hat{\CS}_{n,1-\alpha}=S_6$.
We note that all the quantities required in Procedure 1 can be calculated based on the 
expressions in (\ref{eq:GAN_sample})--(\ref{eq:Qhat_func}) and (\ref{eq:SS_ecdf}).

We also check whether Assumption 3 (required in Theorem 2) seems reasonable -- that is, does the statistic $\sup_{\theta\in\Theta_0}n^{1/2}\hat{Q}_{n}(\theta)$ converge in distribution to some random variable. 
The top-right graph of Figure 1 plots the empirical cdf of 
$\sup_{\theta\in\Theta_0}n^{1/2}\hat{Q}_{n}(\theta)$
based on $200$ separate IID samples of $X_1,\ldots,X_n$ and $Z_1,\ldots,Z_n$
of size $n$, for three different sample sizes 
($n=100$ in blue, $n=1000$ in red, $n=10.000$ in black),
and for the set $\Theta_0$ solved above.
These empirical distributions appear, even for the largest sample size considered, 
non-degenerate (and even continuous), suggesting Assumption 3 holds.

Finally, to empirically verify the result of Theorem 2 in a small Monte Carlo exercise, 
we inspect whether the confidence sets $\hat{\CS}_{n,1-\alpha}$ have correct coverage.
The set-up of our Monte Carlo exercise is as follows.
We consider three different nominal coverage probabilities 
$1-\alpha$ ($80\%$, $90\%$,  $95\%$) and three sample sizes $n$ (100, 500, 1000).
For each sample size $n$, we consider three subsample sizes $b$:
25, 40, and 63 for $n=100$;
77,	144,	and 269 for $n=500$; and
126,  251, 501 for $n=1000$ 
(these correspond to 
$b\approx n^{0.7}$, $b\approx n^{0.8}$, and $b\approx n^{0.9}$
that conform with the requirements $b\to\infty$ and $b/n\to0$ as $n\to\infty$ in Theorem 2).
As for the number of subsamples, we always use $200$ randomly chosen subsamples. 
On each simulation round we form the confidence set $\hat{\CS}_{n,1-\alpha}$ according to Procedure 1
and check whether it covers the entire set $\Theta_0$ or not.
This exercise is repeated $1000$ times and 
Table 1 presents the resulting empirical coverage rates (in $\%$).
Inspection of the results indicates that the choice of subsample size $b$ greatly affects the results; 
this is commonly the case with subsampling, see \citet[Ch 9]{politis1999subsampling}.
It can also be seen that both the sample size $n$ and the subsample size $b$ need to be large enough for the coverage rates to be close to the desired levels.
Nevertheless, for the largest $n$ and $b$ considered, the empirical coverage rates 
$83.4\%$, $92.1\%$, and $95.4\%$
are reasonably close to the desired $80\%$, $90\%$, and $95\%$ nominal levels, 
suggesting that the confidence sets formed using Procedure 1
have the appropriate coverage property stated in Theorem 2.

\begin{figure}[p]

\begin{minipage}[outer sep=0]{\textwidth}
\begin{minipage}[m]{0.2\textwidth}\includegraphics[width=1\textwidth]{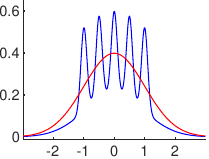}\end{minipage}\hfill
\begin{minipage}[m]{0.3\textwidth}\includegraphics[width=1.9in]{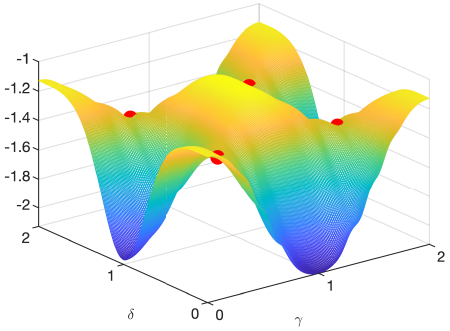}\end{minipage}\hfill
\begin{minipage}[m]{0.3\textwidth}\includegraphics[width=1.9in]{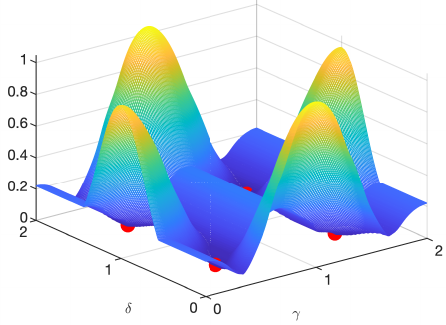}\end{minipage}\hfill
\begin{minipage}[m]{0.2\textwidth}\includegraphics[width=1\textwidth]{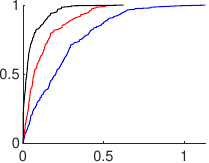}\end{minipage}
\end{minipage}

\vspace*{1pt}

\begin{minipage}[outer sep=0]{\textwidth}
\begin{minipage}[m]{0.245\textwidth}\includegraphics[width=1\textwidth]{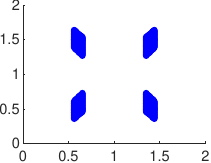}\end{minipage}\hfill
\begin{minipage}[m]{0.245\textwidth}\includegraphics[width=1\textwidth]{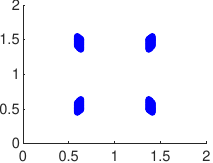}\end{minipage}\hfill
\begin{minipage}[m]{0.245\textwidth}\includegraphics[width=1\textwidth]{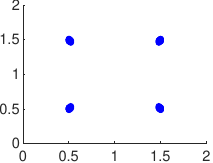}\end{minipage}\hfill
\begin{minipage}[m]{0.245\textwidth}\includegraphics[width=1\textwidth]{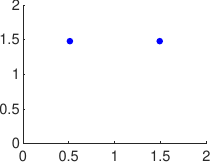}\end{minipage}
\end{minipage}

\vspace*{1pt}

\begin{minipage}[outer sep=0]{\textwidth}
\begin{minipage}[m]{0.245\textwidth}\includegraphics[width=1\textwidth]{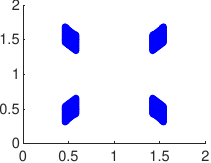}\end{minipage}\hfill
\begin{minipage}[m]{0.245\textwidth}\includegraphics[width=1\textwidth]{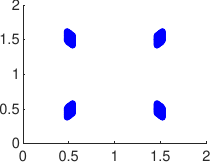}\end{minipage}\hfill
\begin{minipage}[m]{0.245\textwidth}\includegraphics[width=1\textwidth]{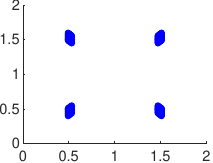}\end{minipage}\hfill
\begin{minipage}[m]{0.245\textwidth}\includegraphics[width=1\textwidth]{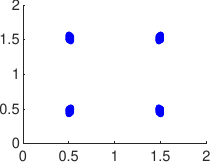}\end{minipage}
\end{minipage}

\vspace*{1pt}

\begin{minipage}[outer sep=0]{\textwidth}
\begin{minipage}[m]{0.245\textwidth}\includegraphics[width=1\textwidth]{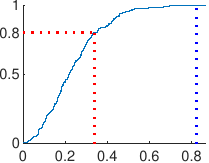}\end{minipage}\hfill
\begin{minipage}[m]{0.245\textwidth}\includegraphics[width=1\textwidth]{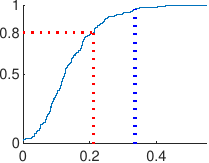}\end{minipage}\hfill
\begin{minipage}[m]{0.245\textwidth}\includegraphics[width=1\textwidth]{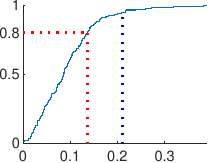}\end{minipage}\hfill
\begin{minipage}[m]{0.245\textwidth}\includegraphics[width=1\textwidth]{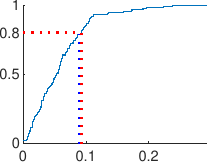}\end{minipage}
\end{minipage}

\vspace*{0pt}
\captionsetup{font=small, width = 0.96\linewidth, labelfont = bf}%
\caption{\label{Fig1}%
Plots for the numerical illustration in Section 5.
\underline{First row:}
(i) Probability density functions of the claw distribution (blue) and the standard normal distribution (red).
(ii)--(iii) Functions $f(\gamma,\delta)$ (second plot) and $Q(\gamma,\delta)$ (third plot), 
with the four solutions in 
$\Theta_0 = \{ (0.5,0.5) , (0.5,1.5) , (1.5,0.5) , (1.5,1.5) \}$
marked with red dots.
(iv) Empirical cdf of $\sup_{\theta\in\Theta_0}n^{1/2}\hat{Q}_{n}(\theta)$
for $n=100$ (blue), $n=1000$ (red), and $n=10.000$ (black).
\underline{Second row:} 
The estimator $\hat{\Theta}_{n}(\tau_n)$ for four different choices of $n$ and $\tau_n$: 
$n=100$ with $\tau_n=0.1n^{-0.25}$ and $\tau_n=0.1n^{-0.49}$ (first two plots) and 
$n=10.000$ with $\tau_n=0.1n^{-0.49}$ and $\tau_n=0$ (third and fourth plots). 
\underline{Third and fourth rows:} 
Illustration of Procedure 1, with sets $S_2$, $S_3$, $S_4$, and the final confidence set $S_6 = \hat{\CS}_{n,1-\alpha}$ on row 3, 
and the corresponding empirical distribution function 
$\hat{L}_{n,b}(S_j,x)$ for $j=2,3,4,6$
(together with the $1-\alpha$ sample quantile $\hat{c}_{n,b,1-\alpha}(S_{j})$, red dotted line, and the quantity $\sup_{\theta\in S_j}n^{1/2}\hat{Q}_{n}(\theta)$, blue dotted line)
on row 4.%
}
\vspace*{20pt}
\end{figure}

\begin{table}[p]
\begin{center}
\renewcommand{\arraystretch}{1.0}
\scalebox{1.0}{
\begin{tabular}{ccccccccccccc}
\hline
\hline\noalign{\medskip}
   &\hspace*{10pt}&& $80\%$  &&\hspace*{10pt}&& $90\%$  &&\hspace*{10pt}&& $95\%$  & \\ 
 	&& \makebox[25pt]{\footnotesize $n\!=\!100$} & \makebox[25pt]{\footnotesize $n\!=\!500$} & \makebox[25pt]{\footnotesize $n\!=\!1000$}   
 	&& \makebox[25pt]{\footnotesize $n\!=\!100$} & \makebox[25pt]{\footnotesize $n\!=\!500$} & \makebox[25pt]{\footnotesize $n\!=\!1000$}   
 	&& \makebox[25pt]{\footnotesize $n\!=\!100$} & \makebox[25pt]{\footnotesize $n\!=\!500$} & \makebox[25pt]{\footnotesize $n\!=\!1000$}   \\
\cline{3-5}\cline{7-9}\cline{11-13}
\noalign{\medskip}
$b\approx n^{0.7}$  && $95.0$  & $97.0$  & $97.2$      && $98.8$  & $98.9$  & $98.9$      && $99.7$  & $99.5$  & $99.6$ \\ 
$b\approx n^{0.8}$  && $90.5$  & $92.4$  & $94.7$      && $96.6$  & $96.1$  & $98.1$      && $98.8$  & $98.5$  & $99.2$ \\ 
$b\approx n^{0.9}$  && $75.1$  & $78.8$  & $83.4$      && $85.9$  & $88.1$  & $92.1$      && $91.3$  & $92.9$  & $95.4$ \\ 
 \noalign{\medskip}
\hline
\hline
\end{tabular}
}
\vspace*{-15pt}
\end{center}
\captionsetup{font=small, width = 0.96\linewidth, labelfont = bf}%
\caption{%
Results of a Monte Carlo study, empirical coverage rate (in $\%$) of the confidence sets $\hat{\CS}_{n,1-\alpha}$. Results are shown for 
three different nominal coverage probabilities ($80\%$, $90\%$,  $95\%$), 
three sample sizes $n$ (100, 500, 1000), and 
three different choices of the subsample size $b$
($b\approx n^{0.7}$, $b\approx n^{0.8}$, and $b\approx n^{0.9}$).%
}
\end{table}

\section{Other minimax problems}

This paper has so far focused on the GAN minimax problem (1) as 
GANs have arguably been one of the key reasons for the recent surge of interest in 
minimax problems in machine learning literature.
In this section we widen our scope and briefly discuss the application of our results to other minimax problems.
We note that the existing literature on general minimax problems is vast, 
originating nearly a century ago with the seminal contribution of 
\citet{vonneumann1928theorie}, 
and it is beyond the scope of this article to review this literature properly;
one could, for instance, see 
\citet{vonneumann2004theory} for minimax problems in game theory, 
and \citet{bertsekas2003convex} for minimax problems in optimization theory.

One standard formulation for general minimax problems appearing 
in the literature is the one used in
\citet{shapiro2008asymptotics}; 
this formulation is also used in the textbooks 
\citet[Sec 5.1.4]{shapiro2009lectures} and \citet[Sec 4.3]{lan2020first}.
The theory developed in the present paper applies, with minor modifications, 
also to this general minimax problem.
Let us outline the required changes. 
Instead of the population and sample GAN problems 
(\ref{eq:GAN_population}) and (\ref{eq:GAN_sample}),
consider the general population and sample minimax problems 
\begin{equation*} \hspace*{-5pt}
\adjustlimits\inf_{\gamma\in\Gamma}\sup_{\delta\in\Delta}f(\gamma,\delta)
\enspace\text{with}\enspace 
f(\gamma,\delta)=\mathbb{E}[F(Z,\gamma,\delta)] 
\quad\text{and}\quad
\adjustlimits\inf_{\gamma\in\Gamma}\sup_{\delta\in\Delta}\hat{f}_{n}(\gamma,\delta)
\enspace\text{with}\enspace
\hat{f}_{n}(\gamma,\delta) = \frac{1}{n}\sum_{i=1}^{n}F(Z_{i},\gamma,\delta),
\end{equation*}
respectively.
No particular interpretation is attached to the 
random vector $Z$ and the parameters $\gamma$ and $\delta$.
Furthermore, replace Assumptions GAN and Smooth GAN with the following two conditions.
\begin{myassumption3*}%
(a) Suppose $Z_{1},\ldots,Z_{n}$ are IID random vectors with the same distribution as $Z$ 
and taking values in $\mathcal{Z}\subseteq\mathbb{R}^{d_{Z}}$.
(b) The set $\Theta=\Gamma\times\Delta\subseteq\mathbb{R}^{d_{\gamma}+d_{\delta}}$ is compact and non-empty.  
(c) The function $F:\mathcal{Z}\times\Gamma\times\Delta\to\mathbb{R}$ is such that 
$F(Z,\gamma,\delta)$ is measurable for all $(\gamma,\delta)\in\Theta$ and continuous on $\Theta$ with probability one.
\end{myassumption3*}
\begin{myassumption4*}%
Assumption Minimax holds, the function $F(z,\theta)$ is continuously differentiable 
on an open convex set $\Theta^{\ast}$ containing $\Theta$ for all $z\in\mathcal{Z}$
with a square integrable derivative
$(\mathbb{E}[\sup_{\theta\in\Theta^{\ast}} |\partial F(Z,\theta)/\partial\theta|^{2}]<\infty)$,
and\/ $\inf_{\theta\in\Theta}\Var[F(Z,\theta)]>0$.
\end{myassumption4*}
With these modifications, all the results of this paper continue to hold also in this general minimax setting. In particular, it is straightforward to see that Theorems 1 and 2 on consistent estimation and confidence sets remain valid. 
In previous literature, 
\citet[Thm 5.9]{shapiro2009lectures} have given a counterpart of Theorem 1(a) in general minimax problems, 
and 
\citet[Thm 3.1]{shapiro2008asymptotics} and \citet[Thm 5.10]{shapiro2009lectures}
have studied the asymptotic distribution of the optimal value in convex-concave minimax problems.
Our results complement these earlier works, in particular by providing confidence sets
for the solutions of the minimax problem. 
A novelty is that all our results hold also in the non-convex and non-concave setting 
(previous literature typically assumes that $f(\cdot,\delta)$ is convex
for all fixed $\delta\in\Delta$ and $f(\gamma,\cdot)$ is concave for
all fixed $\gamma\in\Gamma$)
and in the presence of multiple solutions.

The general minimax formulation considered above is very common and so similar to the GAN problem
that only cosmetic adjustments to our theory were needed.
However, in many examples in the literature the precise formulation of the minimax problem is quite different from the one above. 
For instance in the machine learning literature, 
somewhat different variants of minimax-type problems have been considered at least in 
adversarial learning  
(\citealp{madry2018towards}; \citealp{javanmard2020precise}),
multi-agent reinforcement learning 
(\citealp{wai2018multi}; \citealp{zhang2021multi}),
distributionally robust optimization 
(\citealp{delage2010distributionally}; \citealp{duchi2021learning}),
federated learning (\citealp{mohri2019agnostic}),
signal processing (\citealp{lu2020hybrid}), 
and AUC (area under the ROC curve) maximization (\citealp{ying2016stochastic}).
Adapting the theory of the present paper to such variations of the minimax problem
would require more substantial changes and is left for future research.

\section{Conclusion}

In this paper, we have considered statistical inference for GANs 
and other minimax problems
in the empirically relevant case of multiple solutions. We first considered
the consistency of (approximate) solutions to the sample GAN problem,
and showed that the set of these solutions is a Hausdorff consistent
estimator of the corresponding set of solutions to the population
GAN problem, $\Theta_{0}$. We then presented a subsampling-based iterative
procedure for forming confidence sets for $\Theta_{0}$, and showed
that these confidence sets are conservatively asymptotically consistent.
The consistency result was shown to hold without any restrictions
on the number of solutions to the GAN problem, whereas our results
for confidence sets were given for the common case of multiple but
finite number of solutions. These results extend on the results of
\citet{biau2020some} who considered the case of a single unique solution.
For other general minimax problems, our results 
complement previous works by providing confidence sets
for the solutions of the minimax problem; 
our assumptions allow for the non-convex and non-concave setting
and the presence of multiple solutions.

The present paper considered only the original \citet{goodfellow2014generative}
formulation of the GAN minimax problem 
and the standard formulation of the general minimax problem in Section 6. 
Extensions to other existing GAN variants
(such as the popular Wasserstein GAN of \citealp{arjovsky2017wasserstein})
and to other minimax problems in machine learning and elsewhere 
would be useful. 
Furthermore, the focus in this paper has been theoretical
and exploring the use of our results in practical applications
would be interesting.

\begin{appendix}

\section*{Appendix A:\ Proof of Theorem 1}

\begin{proof}[Proof of Theorem 1]
For completeness, first note that measurability issues and non-emptyness
of $\Theta_{0}$ and $\hat{\Theta}_{n}$ (with probability one) are
discussed in \citet[pp.\ 170--171]{shapiro2009lectures}. For ease
of reference, also note that for any two bounded real-valued functions
$u$ and $v$ defined on a domain $D$ (a subset of some Euclidean
space),
\begin{equation}
| \sup u(x)-\sup v(x) |  \leq  \sup | u(x)-v(x) |,
\;\;\:
| \inf u(x)-\inf v(x) |  \leq  \sup | u(x)-v(x) | , \label{eq:ConsProof_SupInf_Inequalities}
\end{equation}
where all the supremums and infimums are taken over $x\in D$.

(a) To establish the consistency of the optimal value $\hat{V}_{n}$,
it suffices to note that 
\begin{equation}
|\hat{V}_{n}-V_{0}|
=  |\inf_{\gamma\in\Gamma}\hat{\varphi}_{n}(\gamma)-\inf_{\gamma\in\Gamma}\varphi(\gamma)|
\leq  \sup_{\gamma\in\Gamma}|\hat{\varphi}_{n}(\gamma)-\varphi(\gamma)|
\leq  \sup_{\theta\in\Theta}|\hat{f}_{n}(\theta)-f(\theta)|
\overset{p}{\to}0,
\label{eq:ConsProof_VhatN_consistency}
\end{equation}
where we have used the definitions of 
$\hat{V}_{n}$, $V_{0}$, $\hat{\varphi}_{n}(\gamma)$, and $\varphi(\gamma)$, 
the two inequalities in (\ref{eq:ConsProof_SupInf_Inequalities}),
and the uniform convergence in Assumption 1(a). 

We next show that
\begin{equation}
\sup_{\theta\in\Theta}|\hat{Q}_{n}(\theta)-Q(\theta)|\overset{p}{\to}0\quad\text{with the function \ensuremath{Q(\theta)} continuous in \ensuremath{\theta}.}\label{eq:ConsProof_ULLN_Q}
\end{equation}
As $\hat{Q}_{n}(\theta)-Q(\theta)=\hat{\varphi}_{n}(\gamma)-\varphi(\gamma)-(\min\{\hat{f}_{n}(\gamma,\delta),\hat{V}_{n}\}-\min\{f(\gamma,\delta),V_{0}\})$,
the triangle inequality, elementary properties of min and max, and the inequalities in (\ref{eq:ConsProof_VhatN_consistency}) imply that
\begin{align}
|\hat{Q}_{n}(\theta)-Q(\theta)|
&  \leq|\hat{\varphi}_{n}(\gamma)-\varphi(\gamma)|
+\max\{|\hat{f}_{n}(\gamma,\delta)-f(\gamma,\delta)|,|\hat{V}_{n}-V_{0}|\}  \label{eq:ConsProof_QQ_leq_ff} \\
&  \leq2\sup_{\theta\in\Theta}|\hat{f}_{n}(\theta)-f(\theta)| \notag
\end{align}
for all $\theta\in\Theta$. 
Thus by Assumption 1(a), $\sup_{\theta\in\Theta}|\hat{Q}_{n}(\theta)-Q(\theta)|\overset{p}{\to}0$.
As for the continuity in (\ref{eq:ConsProof_ULLN_Q}), by Assumptions
GAN and 1(a) the function $f(\theta)$ is continuous on the compact
set $\Theta$, and thus by Berge's maximum theorem the function $\varphi(\gamma)=\sup_{\delta\in\Delta}f(\gamma,\delta)$
is continuous on $\Gamma$. Consequently, $Q(\theta)$ is continuous.
Also note that the continuity of $Q$ and the definition of $\Theta_{0}$
imply that for all $\epsilon>0$ there exists an $\eta(\epsilon)>0$
such that 
\begin{equation}
\inf_{\theta\in\Theta\setminus\Theta_{0}^{\epsilon}}Q(\theta)\geq\eta(\epsilon),\label{eq:ConsProof_IdCond}
\end{equation}
where $\Theta_{0}^{\epsilon}$ denotes the $\epsilon$-expansion of
the set $\Theta_{0}$ in $\Theta$ defined as $\Theta_{0}^{\epsilon}:=\{\theta\in\Theta:d(\theta,\Theta_{0})\leq\epsilon\}$
and $\Theta\setminus\Theta_{0}^{\epsilon}$ is the complement
of $\Theta_{0}^{\epsilon}$ in $\Theta$. 

Establishing $\sup_{\theta\in\hat{\Theta}_{n}(\tau_{n})}d(\theta,\Theta_{0})\overset{p}{\to}0$
now follows the pattern of a standard consistency proof and relies
on the uniform convergence result (\ref{eq:ConsProof_ULLN_Q}) and
the set-identification condition for $\Theta_{0}$ in (\ref{eq:ConsProof_IdCond}).
To this end, let small $\epsilon_{p},\epsilon_{d}>0$ be arbitrary,
choose an $\eta = \eta(\epsilon_{d})$ as in (\ref{eq:ConsProof_IdCond})
so that $\inf_{\theta\in\Theta\setminus\Theta_{0}^{\epsilon_{d}}}Q(\theta)\geq\eta$
holds, and by (\ref{eq:ConsProof_ULLN_Q}) and Assumption 2(a) choose
an $n_{\epsilon_{p}}$ such that for all $n\geq n_{\epsilon_{p}}$
both $\sup_{\theta\in\Theta} |\hat{Q}_{n}(\theta)-Q(\theta) |\leq\eta/4$
and $\tau_{n}\leq\eta/4$ hold with probability larger than $1-\epsilon_{p}$. We now have 
\begin{align*}
\sup_{\theta\in\hat{\Theta}_{n}(\tau_{n})}Q(\theta)
& \leq\sup_{\theta\in\hat{\Theta}_{n}(\tau_{n})}|\hat{Q}_{n}(\theta)-Q(\theta)|+\sup_{\theta\in\hat{\Theta}_{n}(\tau_{n})}\hat{Q}_{n}(\theta)    \\
& \leq\sup_{\theta\in\Theta} |\hat{Q}_{n}(\theta)-Q(\theta) |+\tau_{n}\leq\eta/2 < \inf_{\theta\in\Theta\setminus\Theta_{0}^{\epsilon_{d}}}Q(\theta),
\end{align*}
and therefore also
$\hat{\Theta}_{n}(\tau_{n})\subseteq\Theta_{0}^{\epsilon_{d}}$
and 
$\sup_{\theta\in\hat{\Theta}_{n}(\tau_{n})}d(\theta,\Theta_{0})\leq\epsilon_{d}$,
for all $n\geq n_{\epsilon_{p}}$ with probability larger than $1-\epsilon_{p}$.
Thus $\sup_{\theta\in\hat{\Theta}_{n}(\tau_{n})}d(\theta,\Theta_{0})\overset{p}{\to}0$.

(b) Assumption 1(b) and (\ref{eq:ConsProof_QQ_leq_ff}) imply that
$\sup_{\theta\in\Theta}n^{1/2} |\hat{Q}_{n}(\theta)-Q(\theta) |=O_{p}(1)$.
Note that
\[
\sup_{\theta\in\Theta_{0}}\hat{Q}_{n}(\theta)\leq\sup_{\theta\in\Theta_{0}}|\hat{Q}_{n}(\theta)-Q(\theta)|+\sup_{\theta\in\Theta_{0}}Q(\theta),
\]
where $\sup_{\theta\in\Theta_{0}}|\hat{Q}_{n}(\theta)-Q(\theta)|\leq\sup_{\theta\in\Theta}|\hat{Q}_{n}(\theta)-Q(\theta)|$
as $\Theta_{0}\subseteq\Theta$ and $\sup_{\theta\in\Theta_{0}}Q(\theta)=0$
by the definition of $\Theta_{0}$. By Assumption 2(b) we have $n^{-1/2}/\tau_{n}\overset{p}{\to}0$ and thus
for any $\epsilon_{p}>0$ we can find an $n_{\epsilon_{p}}$ such
that for all $n\geq n_{\epsilon_{p}}$ 
\[
\sup_{\theta\in\Theta_{0}}\hat{Q}_{n}(\theta)\leq O_{p}(n^{-1/2})=O_{p}(1)(n^{-1/2}/\tau_{n})\tau_{n}\leq\tau_{n}
\]
with probability larger than $1-\epsilon_{p}$. By the definition
of $\hat{\Theta}_{n}(\tau_{n})$ we now have $\Theta_{0}\subseteq\hat{\Theta}_{n}(\tau_{n})$
and thus $\sup_{\theta\in\Theta_{0}}d(\theta,\hat{\Theta}_{n}(\tau_{n}))=0$
(for all $n\geq n_{\epsilon_{p}}$ with probability larger than $1-\epsilon_{p}$). This
shows that $\sup_{\theta\in\Theta_{0}}d(\theta,\hat{\Theta}_{n}(\tau_{n}))\overset{p}{\to}0$.
Combining this with the result in (a), we have established that $d_{H}(\hat{\Theta}_{n}(\tau_{n}),\Theta_{0})\overset{p}{\to}0$. 
\end{proof}

\section*{Appendix B:\ Proofs of Lemma 1 and Theorem 2}

We begin with some preparatory discussion.  
As was outlined in Section 4, 
the statistic\linebreak[4] $\sup_{\theta\in\Theta_{0}}n^{1/2}\hat{Q}_{n}(\theta)$
can be expressed as $n^{1/2}(\phi(\hat{f}_{n})-\phi(f))$ for a suitably
defined map $\phi:l^{\infty}(\Theta)\to\mathbb{R}$ 
(to be given below in the proof of Lemma 1),
this map will be shown to be differentiable in an appropriate sense,
and these facts enable us to apply a suitable functional delta method to obtain
the result of Lemma 1. 
To derive the necessary results, 
it is convenient to define the relevant differentiability concepts for maps between
abstract Banach spaces; the specific spaces used below include $l^{\infty}(\Theta)$,
$C(\Theta)$, and the Euclidean space $\mathbb{R}$. 

Let $X$ and $Y$ be real Banach spaces with norms $\left\Vert \cdot\right\Vert _{X}$
and $\left\Vert \cdot\right\Vert _{Y}$, respectively. Let the domain
$X_{D}\subseteq X$ be some subset of $X$, and consider an arbitrary
map $\phi:X_{D}\to Y$. Let $X_{0}\subseteq X$ be another subset
of $X$. The map $\phi$ is said to be Hadamard directionally differentiable
at $x\in X_{D}$ tangentially to $X_{0}$ with a derivative $\phi_{x}^{\prime}:X_{0}\to Y$
if 
\begin{equation}
\lim_{n\to\infty}\Bigl\Vert\frac{\phi(x+t_{n}h_{n})-\phi(x)}{t_{n}}-\phi_{x}^{\prime}(h)\Bigr\Vert_{Y}=0\label{eq:HDD}
\end{equation}
for all sequences $\{h_{n}\}\subset X$ and $\{t_{n}\}\subset\mathbb{R}_{+}$
such that $t_{n}\downarrow0$ and $h_{n}\to h\in X_{0}$ as $n\to\infty$
and $x+t_{n}h_{n}\in X_{D}$ for all $n$. Two related differentiability
concepts are defined as follows. If the convergence (\ref{eq:HDD})
is required to hold only for all fixed $h_{n}\equiv h\in X_{0}$,
the map $\phi$ is said to be Gateaux directionally differentiable
(at $x\in X_{D}$ tangentially to $X_{0}$). Alternatively, if in
the above definition the requirement ``$\{t_{n}\}\subset\mathbb{R}_{+}$
such that $t_{n}\downarrow0$'' is replaced with ``$\{t_{n}\}\subset\mathbb{R}$
such that $t_{n}\to0$'', the map is said to be (fully) Hadamard
differentiable (at $x\in X_{D}$ tangentially to $X_{0}$). In statistical
literature, these differentiability concepts are discussed and developed
in \citet{reeds1976definition}, \citet{gill1989non}, and Shapiro
\citeyearpar{shapiro1990concepts,shapiro1991asymptotic}, among others. 

Our main interest in these differentiability concepts comes from the
fact that the functional delta method may become applicable. For (fully)
Hadamard differentiable maps this is well known; see,
e.g., \citet[Thm 3.9.4]{vandervaart1996weak}.
For Hadamard directionally differentiable maps the functional
delta method is given in \citet[Thm 2.1]{shapiro1991asymptotic};
recently attention to this result has been drawn for instance by \citet[Thm 2.1]{fang2019inference}
and \citet[Propn 2.1]{carcamo2020directional}. For convenience, we
reproduce this result as the following lemma.
\begin{lem}
Suppose the map $\phi:X_{D}\to Y$ is Hadamard directionally differentiable
at $x\in X_{D}$ tangentially to $X_{0}$ with a derivative $\phi_{x}^{\prime}:X_{0}\to Y$.
Let $x_{1},x_{2},\ldots$ be $X_{D}$-valued random
elements such that $r_{n}(x_{n}-x)\rightsquigarrow x_{0}$
in $X$ for some random element $x_{0}$ taking its values in $X_{0}$
with probability one and for some constants $r_{n}\to\infty$. Then
$r_{n}(\phi(x_{n})-\phi(x))\rightsquigarrow\phi_{x}^{\prime}(x_{0})$
in $Y$. 
\end{lem}
Gateaux directional differentiability is in general too weak of a
differentiability concept to ensure the validity of the functional
delta method. However, when the map $\phi$ is also locally Lipschitz
at $x\in X_{D}$ (in the sense that there exists a $C>0$ such that
for all $x',x''\in X_{D}$ in some neighborhood of $x$, $\lVert\phi(x')-\phi(x'')\rVert_{Y}\leq C\lVert x'-x''\rVert_{X}$),
Gateaux and Hadamard directional differentiability become equivalent
(\citealp[Propn 3.5]{shapiro1990concepts}). This is convenient and
holds for the maps considered in the next lemma. In this lemma, we
let $\Gamma$ and $\Delta$ denote any compact and non-empty Euclidean
sets (which need not equal the sets $\Gamma$ and $\Delta$ elsewhere
in the paper; the same notation is used for convenience). Furthermore, for
any function $x\in C(\Delta)$ we denote $\Delta_{0}(x)=\{\delta_{0}\in\Delta:x(\delta_{0})=\sup_{\delta\in\Delta}x(\delta)\}$,
and for any function $x\in C(\Gamma\times\Delta)$ we let $\bar{x}(\gamma)=\sup_{\delta\in\Delta}x(\gamma,\delta)$
denote the max-function and define the sets 
\begin{align*}
\Theta_{0}(x) & =\{(\gamma_{0},\delta_{0})\in\Gamma\times\Delta:x(\gamma_{0},\delta_{0})=\sup_{\delta\in\Delta}x(\gamma_{0},\delta)=\bar{x}(\gamma_{0})\text{ and }\bar{x}(\gamma_{0})=\inf_{\gamma\in\Gamma}\bar{x}(\gamma)\},\\
\Gamma_{0}(x) & =\{\gamma_{0}\in\Gamma:(\gamma_{0},\delta_{0})\in\Theta_{0}(x)\text{ for some }\delta_{0}\in\Delta\},\vphantom{\sup_{\delta\in\Delta}}\\
\Delta_{0}(\gamma,x) & =\{\delta_{0}\in\Delta:x(\gamma,\delta_{0})=\sup_{\delta\in\Delta}x(\gamma,\delta)\}\quad(\gamma\in\Gamma).
\end{align*}
For any functions $x,y\in l^{\infty}(\Gamma\times\Delta)$, we also allow ourselves
to write 
\[
\underset{(\gamma,\delta)\in\Theta_{0}(x)}{{\textstyle \inf_{\gamma}}\thinspace{\textstyle \sup_{\delta}}}\thinspace y(\gamma,\delta)\quad\text{instead of}\quad\adjustlimits\inf_{\gamma\in\Gamma_{0}(x)}\sup_{\delta\in\Delta_{0}(\gamma,x)}\thinspace y(\gamma,\delta).
\]
Finally, for notational ease the three different maps in parts (a)--(c) of the following lemma
are all denoted by $\phi$ as confusion should not arise.

\begin{lem}
\mbox{}
\begin{lyxlist}{(a)}
\item [{(a)}] The map $\phi:l^{\infty}(\Delta)\to\mathbb{R}$ given by
$\phi(x)=\sup_{\delta\in\Delta}x(\delta)$ for $x\in l^{\infty}(\Delta)$
is Hadamard directionally differentiable at any $x\in C(\Delta)$
tangentially to $C(\Delta)$ with the derivative $\phi'_{x}(h)=\sup_{\delta\in\Delta_{0}(x)}h(\delta)$.\smallskip
\item [{(b)}] The map $\phi:l^{\infty}(\Gamma\times\Delta)\to\mathbb{R}$
given by $\phi(x)=\inf_{\gamma\in\Gamma}\sup_{\delta\in\Delta}x(\gamma,\delta)$
for $x\in l^{\infty}(\Gamma\times\Delta)$ is Hadamard directionally
differentiable at any $x\in C(\Gamma\times\Delta)$ tangentially to
$C(\Gamma\times\Delta)$ with the derivative
\[
\phi'_{x}(h)=\underset{(\gamma,\delta)\in\Theta_{0}(x)}{{\textstyle \inf_{\gamma}}\thinspace{\textstyle \sup_{\delta}}}\thinspace h(\gamma,\delta).
\]
\item [{(c)}] Consider a fixed function $x_{0}\in C(\Gamma\times\Delta)$.
If the set $\Gamma_{0}(x_{0})$ is finite, then the map $\phi:l^{\infty}(\Gamma\times\Delta)\to\mathbb{R}$
given by
\[
\phi(x)=\sup_{(\gamma_{0},\delta_{0})\in\Theta_{0}(x_{0})}\Bigl\{\sup_{\delta\in\Delta}x(\gamma_{0},\delta)-\min\bigl\{ x(\gamma_{0},\delta_{0}),\adjustlimits\inf_{\gamma\in\Gamma}\sup_{\delta\in\Delta}x(\gamma,\delta)\bigr\}\Bigr\}
\]
for $x\in l^{\infty}(\Gamma\times\Delta)$ is Hadamard directionally
differentiable at $x_{0}$ tangentially to $C(\Gamma\times\Delta)$
with the derivative
\[
\phi_{x_{0}}^{\prime}(h)=\sup_{(\gamma_{0},\delta_{0})\in\Theta_{0}(x_{0})}\Bigl\{\sup_{\delta\in\Delta_{0}(\gamma_{0},x_{0})}h(\gamma_{0},\delta)-\min\bigl\{ h(\gamma_{0},\delta_{0}),\underset{(\gamma,\delta)\in\Theta_{0}(x_{0})}{{\textstyle \inf_{\gamma}}\thinspace{\textstyle \sup_{\delta}}}\thinspace h(\gamma,\delta)\bigr\}\Bigr\}.
\]
 
\end{lyxlist}
\end{lem}
The proof of Lemma 3 is given at the end of this appendix.
The result in part (a) is well-known and versions of
it can be found in \citet[Thm 3.1]{shapiro1991asymptotic}, \citet[Lemma S.4.9]{fang2019inference},
and \citet[Corollary 2.3]{carcamo2020directional}, among others.
Result (b) in a convex-concave special case ($x$ and $h$ convex in $\gamma$
and concave in $\delta$) has been obtained by \citet[Propn 2.1]{shapiro2008asymptotics},
but the extension to the present non-convex non-concave case appears to be new
and facilitates our analysis of minimax problems and GANs.
Our main interest is in part (c), the proof of which
relies on parts (a) and (b). The finiteness of $\Gamma_{0}(x_{0})$ is assumed
in part (c) because the mapping taking functions $x$ in $l^{\infty}(\Gamma\times\Delta)$
to functions $\sup_{\delta\in\Delta}x(\gamma,\delta)$, $\gamma\in\Gamma$,
in $l^{\infty}(\Gamma)$ is in general not Hadamard directionally
differentiable without some additional assumptions. 

With the general Lemmas 2 and 3 available, 
we now return to the GAN problem and first give the proof of Lemma 1. 

\begin{proof}[Proof of Lemma 1]
With the set $\Theta_{0}$ as in Section 2,
consider the map $\phi:l^{\infty}(\Gamma\times\Delta)\to\mathbb{R}$
defined, for any function $x\in l^{\infty}(\Gamma\times\Delta)$,
by
\begin{equation*}
\phi(x)=\sup_{(\gamma_{0},\delta_{0})\in\Theta_{0}}\Bigl\{\sup_{\delta\in\Delta}x(\gamma_{0},\delta)-\min\bigl\{ x(\gamma_{0},\delta_{0}),\adjustlimits\inf_{\gamma\in\Gamma}\sup_{\delta\in\Delta}x(\gamma,\delta)\bigr\}\Bigr\} .
\end{equation*}
Noting that the function $\hat{Q}_{n}(\theta)$ in (\ref{eq:Qhat_func}),
here for clarity evaluated at any $\theta_{0}\in\Theta$, can be written as
\begin{equation*}
\hat{Q}_{n}(\theta_{0})=\sup_{\delta\in\Delta}\hat{f}_{n}(\gamma_{0},\delta)-\min\bigl\{\hat{f}_{n}(\gamma_{0},\delta_{0}),\adjustlimits\inf_{\gamma\in\Gamma}\sup_{\delta\in\Delta}\hat{f}_{n}(\gamma,\delta)\bigr\},
\end{equation*}
it is immediate that 
$\sup_{\theta\in\Theta_{0}}n^{1/2}\hat{Q}_{n}(\theta)=n^{1/2}\phi(\hat{f}_{n})$.
Moreover, by the definition of $\Theta_{0}$,\linebreak[4]
$\sup_{\delta\in\Delta}f(\gamma_{0},\delta)=f(\gamma_{0},\delta_{0})=\inf_{\gamma\in\Gamma}\sup_{\delta\in\Delta}f(\gamma,\delta)$
for all $(\gamma_{0},\delta_{0})\in\Theta_{0}$, so that $\phi(f)=0$
and $\sup_{\theta\in\Theta_{0}}n^{1/2}\hat{Q}_{n}(\theta)=n^{1/2}(\phi(\hat{f}_{n})-\phi(f))$.
Applying Lemma 3(c) with the function $x_{0}=f\in C(\Gamma\times\Delta)$,
the sets $\Theta_{0}(f)$ and $\Gamma_{0}(f)$
in Lemma 3(c) reduce to the sets $\Theta_{0}$ and $\Gamma_{0}$ 
in (\ref{eq:Theta0_A}) and in Assumption 5. 
On the other hand, by Assumption 4,
$n^{1/2}(\hat{f}_{n}(\theta)-f(\theta))\rightsquigarrow\mathbb{G}(\theta)$
in $l^{\infty}(\Theta)$ with $\mathbb{G}$ taking values in $C(\Theta)$
with probability one. Combining these facts and using the functional
delta method of Lemma 2 now yields the desired result 
$\sup_{\theta\in\Theta_{0}}n^{1/2}\hat{Q}_{n}(\theta)=n^{1/2}(\phi(\hat{f}_{n})-\phi(f))\rightsquigarrow\phi_{f}^{\prime}(\mathbb{G})$
where
\begin{equation}
\phi_{f}^{\prime}(\mathbb{G})
= \sup_{(\gamma_{0},\delta_{0})\in\Theta_{0}}\Bigl\{  \sup_{\delta\in\Delta_{0}(\gamma_{0},f)}\mathbb{G}(\gamma_{0},\delta)
- \min\bigl\{\mathbb{G}(\gamma_{0},\delta_{0}),\underset{(\gamma,\delta)\in\Theta_{0}}{{\textstyle \inf_{\gamma}}\thinspace{\textstyle \sup_{\delta}}}\thinspace\mathbb{G}(\gamma,\delta)\bigr\}\Bigr\}.\label{eq:LimitingRV}
\end{equation}
\end{proof}

Expression (\ref{eq:LimitingRV}) gives the limiting distribution $L$ in Assumption 3.
In general this distribution is rather complicated.
To illustrate, consider the case with multiple but finite number of
distinct solutions mentioned after Theorem 2 ($\Theta_{0}=\{(\gamma_{01},\delta_{01}),\ldots,(\gamma_{0K},\delta_{0K})\}$
with $K>1$ and $\gamma_{0i}\neq\gamma_{0j}$, $\delta_{0i}\neq\delta_{0j}$
for all $i\neq j$). It may be helpful to note that (\ref{eq:LimitingRV})
can equivalently be expressed as 
\begin{align*}
\adjustlimits\sup_{\gamma_{0}\in\Gamma_{0}}\sup_{\delta_{0}\in\Delta_{0}(\gamma_{0},f)}
\Bigl\{
\max\bigl\{ 
& \sup_{\delta\in\Delta_{0}(\gamma_{0},f)}\mathbb{G}(\gamma_{0},\delta)-\mathbb{G}(\gamma_{0},\delta_{0}), \\
&\qquad \sup_{\delta\in\Delta_{0}(\gamma_{0},f)}\mathbb{G}(\gamma_{0},\delta)-\adjustlimits\inf_{\gamma\in\Gamma_{0}}\sup_{\delta\in\Delta_{0}(\gamma,f)}\thinspace\mathbb{G}(\gamma,\delta)\bigr\}
\Bigr\}
\end{align*}
and in the said special case this expression reduces to 
\begin{equation}
\max_{k=1,\ldots,K}\mathbb{G}(\gamma_{0k},\delta_{0k})-\min_{k=1,\ldots,K}\mathbb{G}(\gamma_{0k},\delta_{0k})\label{eq:LimitingRV_SpecialCase}
\end{equation}
which also satisfies the continuity requirement of Theorem 2 
due to \citet[Thm 11.1]{davydov1998local}
(as the mapping from 
$(\mathbb{G}(\gamma_{01},\delta_{01}),\ldots, \linebreak[1] \mathbb{G}(\gamma_{0K},\delta_{0K}))$
to (\ref{eq:LimitingRV_SpecialCase}) is convex and 
$\inf_{\theta\in\Theta_{0}}\mathbb{E}[(\mathbb{G}(\theta))^{2}] \linebreak[4] >0$ by Assumption 4).

When the set $\Theta_{0}$ is a singleton, 
the derivative map $\phi_{f}^{\prime}$ in (\ref{eq:LimitingRV}) is identically zero 
and the corresponding limiting distribution $L$ is degenerate.
To obtain an informative limiting distribution in this case, 
our approach would need to be modified.
One possibility is to apply a higher-order functional delta method 
based on second-order directional derivatives 
similarly as in the second-order analysis for minimization problems 
in \citet[Sec 5.1.3]{shapiro2009lectures}.
However, we do not pursue this further as 
in the case of a single solution confidence sets can be directly formed based on the 
asymptotic normality result of \citet[Thm 4.3]{biau2020some}.

We are now ready to prove Theorem 2.

\begin{proof}[Proof of Theorem 2]
Suppose the confidence sets $\hat{\CS}_{n,1-\alpha}$ ($n=1,2,\ldots$) are formed
using Procedure 1 so that $\hat{\CS}_{n,1-\alpha}=S_{\hat{j}_{n}}$
for some indices $\hat{j}_{n}$. 
First consider the probabilities $\mathbb{P}[\Theta_{0}\nsubseteq\hat{\CS}_{n,1-\alpha}]$.
In the event $\Theta_{0}\nsubseteq\hat{\CS}_{n,1-\alpha}$,
necessarily $\hat{j}_{n}\geq2$ (as $\Theta_{0}\subseteq\Theta=S_{1}$)
and, by the construction of $\hat{\CS}_{n,1-\alpha}$ in Procedure 1, 
there exists an index $\hat{i}_{n}\in\{2,\ldots,\hat{j}_{n}\}$ such that 
$\Theta_{0}\subseteq{S\/}_{\hat{i}_{n}-1}$ but $\Theta_{0}\nsubseteq{S\/}_{\hat{i}_{n}}$.
For such an index  $\hat{i}_{n}$, there exists a $\theta^{\bullet}$ such that $\theta^{\bullet}\in\Theta_{0}$
but $\theta^{\bullet}\notin{S\/}_{\hat{i}_{n}}$ and then 
\[
\sup_{\theta\in\Theta_{0}}n^{1/2}\hat{Q}_{n}(\theta)\geq n^{1/2}\hat{Q}_{n}(\theta^{\bullet})>\hat{c}_{n,b,1-\alpha}({S\/}_{\hat{i}_{n}-1})\geq\hat{c}_{n,b,1-\alpha}(\Theta_{0}),
\]
where the first inequality holds as $\theta^{\bullet}\in\Theta_{0}$,
the second inequality is due to the definition of ${S\/}_{\hat{i}_{n}}$,
and the third inequality follows from the fact that $\Theta_{0}\subseteq{S\/}_{\hat{i}_{n}-1}$
(this last fact implies that $\hat{L}_{n,b}(\Theta_{0},x)\geq\hat{L}_{n,b}({S\/}_{\hat{i}_{n}-1},x)$
for all $x$, and thus the third inequality in the previous display
follows). Therefore it holds that 
\begin{equation}
\mathbb{P}[\Theta_{0}\subseteq\hat{\CS}_{n,1-\alpha}]
= 1 - \mathbb{P}[\Theta_{0}\nsubseteq\hat{\CS}_{n,1-\alpha}]
\geq \mathbb{P}[ \sup_{\theta\in\Theta_{0}}n^{1/2}\hat{Q}_{n}(\theta)
\leq \hat{c}_{n,b,1-\alpha}(\Theta_{0}) ].
\label{eq:Thm2proof1}
\end{equation}

Now consider the sample quantiles $\hat{c}_{n,b,1-\alpha}(\Theta_{0})$ 
appearing in (\ref{eq:Thm2proof1}) and the associated empirical distribution function 
$\hat{L}_{n,b}(\Theta_{0},x)$.
By Assumption 3 and the conditions imposed on the subsampling size
$b$ in Theorem 2, the requirements of Theorem 2.2.1(i) of \citet{politis1999subsampling}
are satisfied so that $\hat{L}_{n,b}(\Theta_{0},x)\overset{p}{\to}L(x)$
for any $x$ that is a continuity point of the cdf
$L(\cdot)$. In particular, for any $\epsilon>0$ for which $c_{1-\alpha}-\epsilon$
is such a continuity point, $\hat{L}_{n,b}(\Theta_{0},c_{1-\alpha}-\epsilon)\overset{p}{\to}L(c_{1-\alpha}-\epsilon)<1-\alpha$
where the inequality holds by the definition of $c_{1-\alpha}$. Together
with the definition of $\hat{c}_{n,b,1-\alpha}(\Theta_{0})$, this
implies that $\lim_{n\to\infty}\mathbb{P}[c_{1-\alpha}-\epsilon\leq\hat{c}_{n,b,1-\alpha}(\Theta_{0})]=1$.
Elementary probability rules now imply that
\begin{align*}
 & \mathbb{P}[\sup_{\theta\in\Theta_{0}}n^{1/2}\hat{Q}_{n}(\theta)\leq\hat{c}_{n,b,1-\alpha}(\Theta_{0})]\\
 & \quad \geq\mathbb{P}[\sup_{\theta\in\Theta_{0}}n^{1/2}\hat{Q}_{n}(\theta)\leq\hat{c}_{n,b,1-\alpha}(\Theta_{0})\text{ and }c_{1-\alpha}-\epsilon\leq\hat{c}_{n,b,1-\alpha}(\Theta_{0})]\\
 & \quad\geq\mathbb{P}[\sup_{\theta\in\Theta_{0}}n^{1/2}\hat{Q}_{n}(\theta)\leq c_{1-\alpha}-\epsilon\text{ and }c_{1-\alpha}-\epsilon\leq\hat{c}_{n,b,1-\alpha}(\Theta_{0})]\\
 & \quad=\mathbb{P}[\sup_{\theta\in\Theta_{0}}n^{1/2}\hat{Q}_{n}(\theta)\leq c_{1-\alpha}-\epsilon]+o(1),
\end{align*}
where the $o(1)$ term converges to zero as $n\to\infty$ due to $\lim_{n\to\infty}\mathbb{P}[c_{1-\alpha}-\epsilon\leq\hat{c}_{n,b,1-\alpha}(\Theta_{0})]=1$.
Therefore
\begin{equation}
\liminf_{n\to\infty}\mathbb{P}[\sup_{\theta\in\Theta_{0}}n^{1/2}\hat{Q}_{n}(\theta)\leq\hat{c}_{n,b,1-\alpha}(\Theta_{0})]\geq L(c_{1-\alpha}-\epsilon).
\label{eq:Thm2proof2}
\end{equation}
Using (\ref{eq:Thm2proof1}) and (\ref{eq:Thm2proof2}) and 
choosing a sequence $\epsilon_{k}>0$, $k=1,2,\ldots$, such that
$c_{1-\alpha}-\epsilon_{k}$ is a continuity point of $L(\cdot)$
for each $k=1,2,\ldots$ and such that $\epsilon_{k}\to0$ as $k\to\infty$
now yields the result of Theorem 2.
\end{proof}

We close this appendix with the proof of Lemma 3.

\begin{proof}[Proof of Lemma 3]
We omit the proof of (the well-known) part (a) and focus on the new parts (b) and (c).
To prove part (b), first note that for any $x',x'' \in l^{\infty}(\Gamma\times\Delta)$ we have
(making use of the two inequalities in (\ref{eq:ConsProof_SupInf_Inequalities})) that
\begin{equation}
|\phi(x')-\phi(x'')| 
\leq \sup_{\gamma\in\Gamma}|\sup_{\delta\in\Delta}x'(\gamma,\delta)  
-  \sup_{\delta\in\Delta}x''(\gamma,\delta)| 
\leq \sup_{(\gamma,\delta)\in\Gamma\times\Delta}|x'(\gamma,\delta)  -  x''(\gamma,\delta)|\label{eq:HDD_infsup_Lipschitz}
\end{equation}
so that the map $\phi$ is Lipschitz and it therefore suffices to
consider Gateaux directional differentiability (recall the remarks following Lemma 2). 
To this end, pick any
$x,h\in C(\Gamma\times\Delta)$ and a sequence $\{t_{n}\}\subset\mathbb{R}_{+}$
such that $t_{n}\downarrow0$ as $n\to\infty$. Denote $x_{(n)}=x+t_{n}h$.
For any $\gamma$ and each $n$, set $\Delta_{0}(\gamma,x)=\{\delta_{0}\in\Delta:x(\gamma,\delta_{0})=\sup_{\delta\in\Delta}x(\gamma,\delta)\}$
and $\Delta_{0}(\gamma,x_{(n)})=\{\delta_{0}\in\Delta:x_{(n)}(\gamma,\delta_{0})=\sup_{\delta\in\Delta}x_{(n)}(\gamma,\delta)\}$;
by continuity and compactness, these sets are non-empty. 
In what follows, we consider the difference quotient 
\[
D_{n}=t_{n}^{-1}\Bigl[\thinspace\adjustlimits\inf_{\gamma\in\Gamma}\sup_{\delta\in\Delta}x_{(n)}(\gamma,\delta)-\adjustlimits\inf_{\gamma\in\Gamma}\sup_{\delta\in\Delta}x(\gamma,\delta)\thinspace\Bigr]
\]
and prove that 
\begin{equation}
\limsup_{n\to\infty}D_{n}
\leq
\underset{(\gamma,\delta)\in\Theta_{0}(x)}{{\textstyle \inf_{\gamma}}\thinspace{\textstyle \sup_{\delta}}}\thinspace h(\gamma,\delta)
\leq\liminf_{n\to\infty}D_{n},
\label{eq:HDD_infsup_proof1}
\end{equation}
implying the desired result.

To study the limes superior of $D_{n}$, first choose an arbitrary
$\gamma^{*}\in\Gamma$ such that $(\gamma^{*},\delta)\in\Theta_{0}(x)$
for some $\delta\in\Delta$. Then, for each $n$, choose a $\delta_{n}\in\Delta_{0}(\gamma^{*},x_{(n)})$
(we suppress the dependence of $\delta_{n}$ on $\gamma^{*}$). By
the definitions of the infimum and $\Delta_{0}(\gamma^{*},x_{(n)})$,
\[
\adjustlimits\inf_{\gamma\in\Gamma}\sup_{\delta\in\Delta}x_{(n)}(\gamma,\delta)\leq\sup_{\delta\in\Delta}x_{(n)}(\gamma^{*},\delta)=x_{(n)}(\gamma^{*},\delta_{n})
\]
and by the definitions of $\gamma^{*}$, $\Theta_{0}(x)$, and the
supremum, 
\[
\adjustlimits\inf_{\gamma\in\Gamma}\sup_{\delta\in\Delta}x(\gamma,\delta)=\sup_{\delta\in\Delta}x(\gamma^{*},\delta)\geq x(\gamma^{*},\delta_{n}).
\]
Therefore 
\[
\adjustlimits\inf_{\gamma\in\Gamma}\sup_{\delta\in\Delta}x_{(n)}(\gamma,\delta)-\adjustlimits\inf_{\gamma\in\Gamma}\sup_{\delta\in\Delta}x(\gamma,\delta)\leq x_{(n)}(\gamma^{*},\delta_{n})-x(\gamma^{*},\delta_{n})=t_{n}h(\gamma^{*},\delta_{n})
\]
so that $D_{n}\leq h(\gamma^{*},\delta_{n})$.
Now let $\delta_{n_{k}}$be an arbitrary convergent subsequence of $\delta_{n}$ 
(compactness of $\Delta$ ensures the existence of such subsequences) 
and suppose $\delta_{n_{k}}\to\delta^{*}$ with $\delta^{*}\in\Delta$ as $k\to\infty$. 
As (for all $k$ and all $\delta\in\Delta$) 
\[
x(\gamma^{*},\delta_{n_{k}})+t_{n_{k}}h(\gamma^{*},\delta_{n_{k}})
=x_{(n_{k})}(\gamma^{*},\delta_{n_{k}})
\geq x_{(n_{k})}(\gamma^{*},\delta)  =  x(\gamma^{*},\delta)+t_{n_{k}}h(\delta),
\]
it also holds (for all $k$ and all $\delta\in\Delta$) that
\[
x(\gamma^{*},\delta^{*})
\geq 
x(\gamma^{*},\delta)+t_{n_{k}}h(\gamma^{*},\delta)-t_{n_{k}}h(\gamma^{*},\delta_{n_{k}})+x(\gamma^{*},\delta^{*})-x(\gamma^{*},\delta_{n_{k}}).
\]
Letting $k$ tend to infinity, continuity of $x$ implies that 
$x(\gamma^{*},\delta^{*})\geq x(\gamma^{*},\delta)$
for all $\delta\in\Delta$. Therefore, by the definition of $\gamma^{*}$,
\[
x(\gamma^{*},\delta^{*})=\sup_{\delta\in\Delta}x(\gamma^{*},\delta)=\adjustlimits\inf_{\gamma\in\Gamma}\sup_{\delta\in\Delta}x(\gamma,\delta)
\]
so that $\delta^{*}$ is such that $\delta^{*}\in\Delta_{0}(\gamma^{*},x)$
and $(\gamma^{*},\delta^{*})\in\Theta_{0}(x)$. 
Continuity of $h$ now implies that $h(\gamma^{*},\delta_{n_{k}})\to h(\gamma^{*},\delta^{*})\leq\sup_{\delta\in\Delta_{0}(\gamma^{*},x)}h(\gamma^{*},\delta)$
as $k\to\infty$. As the convergent subsequence $\delta_{n_{k}}$was
arbitrary, $\limsup_{n\to\infty}h(\gamma^{*},\delta_{n})\leq\sup_{\delta\in\Delta_{0}(\gamma^{*},x)}h(\gamma^{*},\delta)$.
Therefore also $\limsup_{n\to\infty}D_{n} \linebreak[4] \leq  \sup_{\delta\in\Delta_{0}(\gamma^{*},x)}h(\gamma^{*},\delta)$.
As $\gamma^{*}\in\Gamma$ was arbitrary except for satisfying $(\gamma^{*},\delta)\in\Theta_{0}(x)$
for some $\delta\in\Delta$, 
\[
\limsup_{n\to\infty}D_{n}\leq\underset{(\gamma,\delta)\in\Theta_{0}(x)}{{\textstyle \inf_{\gamma}}\thinspace{\textstyle \sup_{\delta}}}\thinspace h(\gamma,\delta).
\]

For the limes inferior of $D_{n}$, it suffices to note that for any $\gamma\in\Gamma$
and for all $n$
\[
\sup_{\delta\in\Delta}x_{(n)}(\gamma,\delta)\geq\sup_{\delta\in\Delta_{0}(\gamma,x)}x_{(n)}(\gamma,\delta)=\sup_{\delta\in\Delta}x(\gamma,\delta)+t_{n}\sup_{\delta\in\Delta_{0}(\gamma,x)}h(\gamma,\delta)
\]
so that (for any $\gamma\in\Gamma$ and for all $n$)
\begin{align*}
&t_{n}^{-1}[\sup_{\delta\in\Delta}x_{(n)}(\gamma,\delta)-\adjustlimits\inf_{\gamma\in\Gamma}\sup_{\delta\in\Delta}x(\gamma,\delta)] 	\\
&\qquad \geq t_{n}^{-1}[\sup_{\delta\in\Delta}x(\gamma,\delta)-\adjustlimits\inf_{\gamma\in\Gamma}\sup_{\delta\in\Delta}x(\gamma,\delta)]
+ \!\sup_{\delta\in\Delta_{0}(\gamma,x)}\!h(\gamma,\delta)\\
&\qquad \geq t_{n}^{-1}  [  \sup_{\delta\in\Delta}x(\gamma,\delta)-\adjustlimits\inf_{\gamma\in\Gamma}\sup_{\delta\in\Delta}x(\gamma,\delta)  ]	\\
&\qquad\qquad\qquad+  
  \sup_{\delta\in\Delta_{0}(\gamma,x)}h(\gamma,\delta)-\underset{(\gamma,\delta)\in\Theta_{0}(x)}{{\textstyle \inf_{\gamma}}\thinspace{\textstyle \sup_{\delta}}}\thinspace h(\gamma,\delta)   
+\underset{(\gamma,\delta)\in\Theta_{0}(x)}{{\textstyle \inf_{\gamma}}\thinspace{\textstyle \sup_{\delta}}}\thinspace h(\gamma,\delta)\\
&\qquad \geq\underset{(\gamma,\delta)\in\Theta_{0}(x)}{{\textstyle \inf_{\gamma}}\thinspace{\textstyle \sup_{\delta}}}\thinspace h(\gamma,\delta).
\end{align*}
Taking an infimum over $\gamma\in\Gamma$ and then the limes inferior as $n\to\infty$ 
establishes the latter inequality in (\ref{eq:HDD_infsup_proof1}), hence  
completing the proof of part (b).

To prove part (c), we first use the former inequality in 
(\ref{eq:ConsProof_SupInf_Inequalities}) and the triangle
inequality to note that for any $x',x''\in l^{\infty}(\Gamma\times\Delta)$
\begin{align*}
& |\phi(x')-\phi(x'')| \\*
& = \;\Bigr|\sup_{(\gamma_{0},\delta_{0})\in\Theta_{0}(x_{0})}\Bigl\{\sup_{\delta\in\Delta}x'(\gamma_{0},\delta)-\min\bigl\{ x'(\gamma_{0},\delta_{0}),\adjustlimits\inf_{\gamma\in\Gamma}\sup_{\delta\in\Delta}x'(\gamma,\delta)\bigr\}\Bigr\}\\*
 & \quad -\sup_{(\gamma_{0},\delta_{0})\in\Theta_{0}(x_{0})}\Bigl\{\sup_{\delta\in\Delta}x''(\gamma_{0},\delta)-\min\bigl\{ x''(\gamma_{0},\delta_{0}),\adjustlimits\inf_{\gamma\in\Gamma}\sup_{\delta\in\Delta}x''(\gamma,\delta)\bigr\}\Bigr\}\Bigl|\\
 & \leq\sup_{(\gamma_{0},\delta_{0})\in\Theta_{0}(x_{0})}\Bigr|\sup_{\delta\in\Delta}x'(\gamma_{0},\delta)-\sup_{\delta\in\Delta}x''(\gamma_{0},\delta)\Bigl|\\*
 & \quad+\sup_{(\gamma_{0},\delta_{0})\in\Theta_{0}(x_{0})}\Bigr|\min\bigl\{ x'(\gamma_{0},\delta_{0}),\adjustlimits\inf_{\gamma\in\Gamma}\sup_{\delta\in\Delta}x'(\gamma,\delta)\bigr\}-\min\bigl\{ x''(\gamma_{0},\delta_{0}),\adjustlimits\inf_{\gamma\in\Gamma}\sup_{\delta\in\Delta}x''(\gamma,\delta)\bigr\}\Bigl|.
\end{align*}
Using the former inequality in (\ref{eq:ConsProof_SupInf_Inequalities})
again it is easy to see that the first term on the majorant side of
the previous display is dominated by $\sup_{(\gamma,\delta)\in\Gamma\times\Delta}|x'(\gamma,\delta)-x''(\gamma,\delta)|$.
Also the second term is dominated by the same quantity, as can be
seen by using the elementary inequality $|\min\{a,b\}-\min\{c,d\}|\leq\max\{|a-c|,|b-d|\}$
and the inequalities in (\ref{eq:HDD_infsup_Lipschitz}). Thus the
map $\phi$ in part (c) is Lipschitz and it suffices to consider Gateaux
directional differentiability. 

Pick any $h\in C(\Gamma\times\Delta)$, a sequence $\{t_{n}\}\subset\mathbb{R}_{+}$
such that $t_{n}\downarrow0$ as $n\to\infty$, and denote $x_{(n)}=x_{0}+t_{n}h$.
Note that by the definition of $\Theta_{0}(x_{0})$, 
\[
\sup_{\delta\in\Delta}x_{0}(\gamma_{0},\delta)=x_{0}(\gamma_{0},\delta_{0})=\adjustlimits\inf_{\gamma\in\Gamma}\sup_{\delta\in\Delta}x_{0}(\gamma,\delta)\quad\text{for all }(\gamma_{0},\delta_{0})\in\Theta_{0}(x_{0})
\]
so that $\phi(x_{0})=0$. These equalities also imply that $\phi(x_{0}+t_{n}h)=\phi(x_{(n)})$
can be written as
\begin{align*}
&  \sup_{(\gamma_{0},\delta_{0})\in\Theta_{0}(x_{0})}
\Bigl\{
\sup_{\delta\in\Delta}x_{(n)}(\gamma_{0},\delta)-\sup_{\delta\in\Delta}x_{0}(\gamma_{0},\delta)      \\*
& \qquad\qquad\qquad\;\; -\min\bigl\{    x_{(n)}(\gamma_{0},\delta_{0})-x_{0}(\gamma_{0},\delta_{0}),\adjustlimits\inf_{\gamma\in\Gamma}\sup_{\delta\in\Delta}x_{(n)}(\gamma,\delta)-\adjustlimits\inf_{\gamma\in\Gamma}\sup_{\delta\in\Delta}x_{0}(\gamma,\delta)    \bigr\}
\Bigr\}
\end{align*}
and therefore the difference quotient $t_{n}^{-1}\{\phi(x_{0}+t_{n}h)-\phi(x_{0})\}=t_{n}^{-1}\{\phi(x_{(n)})-\phi(x_{0})\}$
can be written as
\begin{align*}
&  \sup_{(\gamma_{0},\delta_{0})\in\Theta_{0}(x_{0})}
\Bigl\{ 
t_{n}^{-1} \bigl\{  \sup_{\delta\in\Delta}x_{(n)}(\gamma_{0},\delta)
	-\sup_{\delta\in\Delta}x_{0}(\gamma_{0},\delta)  \bigr\}	\\
& \qquad\qquad\qquad  -\min\bigl\{ h(\gamma_{0},\delta_{0}),t_{n}^{-1}\bigl\{\adjustlimits\inf_{\gamma\in\Gamma}\sup_{\delta\in\Delta}x_{(n)}(\gamma,\delta)-\adjustlimits\inf_{\gamma\in\Gamma}\sup_{\delta\in\Delta}x_{0}(\gamma,\delta)\bigr\}\bigr\}
\Bigr\}.
\end{align*}
Introducing the shorthand notation 
\begin{align*}
A_{n}(\gamma) & =t_{n}^{-1}\bigl\{\sup_{\delta\in\Delta}x_{(n)}(\gamma,\delta)-\sup_{\delta\in\Delta}x_{0}(\gamma,\delta)\bigr\},\\ 
B_{n} &=t_{n}^{-1}\bigl\{\adjustlimits\inf_{\gamma\in\Gamma}\sup_{\delta\in\Delta}x_{(n)}(\gamma,\delta)-\adjustlimits\inf_{\gamma\in\Gamma}\sup_{\delta\in\Delta}x_{0}(\gamma,\delta)\bigr\},\\
C_{n}(\gamma) & =\sup_{\delta\in\Delta_{0}(\gamma,x_{0})}h(\gamma,\delta),
\qquad \quad
D_{n}=\underset{(\gamma,\delta)\in\Theta_{0}(x_{0})}{{\textstyle \inf_{\gamma}}\thinspace{\textstyle \sup_{\delta}}}\thinspace h(\gamma,\delta),
\end{align*}
and using the inequality in (\ref{eq:ConsProof_SupInf_Inequalities})
and the triangle inequality allow us to write 
\begin{align*}
&  \bigl|t_{n}^{-1}\bigl\{\phi(x_{(n)})-\phi(x_{0})\bigr\}-\phi_{x_{0}}^{\prime}(h)\bigr|	\\*
& \qquad=\bigl|\sup_{(\gamma_{0},\delta_{0})\in\Theta_{0}(x_{0})}\bigl\{ A_{n}(\gamma_{0})-\min\{h(\gamma_{0},\delta_{0}),B_{n}\}\bigr\} \\*
&\qquad\qquad\qquad -\sup_{(\gamma_{0},\delta_{0})\in\Theta_{0}(x_{0})}\bigl\{ C_{n}(\gamma_{0})-\min\{h(\gamma_{0},\delta_{0}),D_{n}\}\bigr\}\bigr|	\\
& \qquad\leq\sup_{(\gamma_{0},\delta_{0})\in\Theta_{0}(x_{0})}\bigl|A_{n}(\gamma_{0})-C_{n}(\gamma_{0})\bigr| \\
&\qquad\qquad\qquad +\sup_{(\gamma_{0},\delta_{0})\in\Theta_{0}(x_{0})}\bigl|\min\{h(\gamma_{0},\delta_{0}),B_{n}\}-\min\{h(\gamma_{0},\delta_{0}),D_{n}\}\bigr| .
\end{align*}
By the definition of the set $\Gamma_{0}(x_{0})$, the former term on the majorant side of the previous display can be written as
$\sup_{\gamma_{0}\in\Gamma_{0}(x_{0})} | A_{n}(\gamma_{0})-C_{n}(\gamma_{0}) |$.
For each fixed $\gamma_{0}\in\Gamma_{0}(x_{0})$ at a time, 
$| A_{n}(\gamma_{0})-C_{n}(\gamma_{0}) | \to 0$
as $n\to\infty$ by the definitions of $A_{n}(\gamma)$ and $C_{n}(\gamma)$,
part (a) of this lemma, and the definition of Hadamard directional
differentiability. By assumption, the set $\Gamma_{0}(x_{0})$ is
finite, implying that $\sup_{\gamma_{0}\in\Gamma_{0}(x_{0})} | A_{n}(\gamma_{0})-C_{n}(\gamma_{0}) | \to 0$
as $n\to\infty$ (note that this result does not generally hold if
$\Gamma_{0}(x_{0})$ is infinite). On the other hand, the elementary
inequality $|\min\{a,b\}-\min\{c,d\}|\leq\max\{|a-c|,|b-d|\}$ can
be used to show that the latter term on the majorant side of the previous
display is dominated by $|B_{n}-D_{n}|$. By the definitions of $B_{n}$
and $D_{n}$, part (b) of this lemma, and the definition of Hadamard
directional differentiability, $|B_{n}-D_{n}|\to0$ as $n\to\infty$.
Overall, we have shown that $|t_{n}^{-1}\{\phi(x_{(n)})-\phi(x_{0})\}-\phi_{x_{0}}^{\prime}(h)| \to 0$
as $n\to\infty$, which completes the proof of part (c).
\end{proof}

\section*{Appendix C:\ Sufficient conditions for Assumptions 1 and 4}

Assumption GAN and condition $\mathbb{E}[\sup_{\theta\in\Theta}|F(X,Z,\theta)|]<\infty$
imply Assumption 1(a) by \citet[Sec 9.17]{hoffmannjorgensen1994probabilityII}
(cf.~also the more transparent Example 19.8 and Theorem 19.4 of \citet{vanderVaart1998asymptotic}
where slightly stronger assumptions are made). 

When Assumption Smooth GAN holds, the
mean value theorem implies that the Lipschitz condition 
$| F(x,z,\theta_{1})-F(x,z,\theta_{2}) | 
\leq 
\sup_{\theta\in\Theta^{\ast}} 
| \partial F(x,z,\theta)/\partial\theta | \, | \theta_{1}-\theta_{2} |$
holds for all $\theta_{1},\theta_{2}\in\Theta$
where 
$\mathbb{E}[
\sup_{\theta\in\Theta^{\ast}}
|\partial F(X,Z,\theta)/\partial\theta|^{2}]<\infty$
by assumption. Example 19.7 and Theorem 19.5 of \citet{vanderVaart1998asymptotic}
now imply the weak convergence in Assumption 4; 
moreover, it holds that $\mathbb{E}[\sup_{\theta\in\Theta}(F(X,Z,\theta))^{2}] \linebreak[3]<\infty$
(\citealp[p.~129]{vandervaart1996weak}).
Condition \linebreak[4] $\inf_{\theta\in\Theta_{0}}\mathbb{E}[(\mathbb{G}(\theta))^{2}]>0$
follows from the requirement $\inf_{\theta\in\Theta}\Var[F(X,Z,\theta)]>0$
in Assumption Smooth GAN. 
Therefore Assumptions 1(b) and 4 hold.

\end{appendix}

\setstretch{1.04}

\bibliographystyle{elsarticle-harv} 
\bibliography{GAN,GAN_v2_additions}

\end{document}